\documentclass{article}
\usepackage{arxiv}
\usepackage{amsmath}
\usepackage[utf8]{inputenc} 
\usepackage[T1]{fontenc}    
\usepackage{hyperref}       
\usepackage{url}            
\usepackage{booktabs}       
\usepackage{adjustbox}  
\usepackage{titlesec}
\usepackage{amsfonts}       
\usepackage{nicefrac}       
\usepackage{microtype}      
\usepackage{graphicx}
\usepackage{float}
\usepackage{afterpage}
\usepackage{booktabs, calc}
\usepackage{filecontents}
\usepackage{multicol}
\usepackage{multirow}
\usepackage{indentfirst}
\usepackage{algorithm}
\usepackage{algpseudocode}
\usepackage{subcaption}
\usepackage{tikz}
\usetikzlibrary{positioning, fit, shapes.geometric, arrows.meta, calc, backgrounds}
\newtheorem{example}{Example}
\DeclareMathAlphabet\mathbfcal{OMS}{cmsy}{b}{n} 
\usepackage{lineno}

\title{MATLAB-Based Layerwise Self-Adaptive Physics-Informed Neural Network in Applications to Multidimensional Coupled Burgers' Equations with High Reynolds Numbers}

\author{
Harish Bhatt \\
  Department of Mathematics\\
  Utah Valley University\\
 Orem, UT 84058 \\
  \texttt{harish.bhatt@uvu.edu} \\ 
  \And
  Xi Chen  \\
Department of Computer Science\\
Utah Valley University\\
Orem, UT 84058\\
  \texttt{Xi.Chen@uvu.edu} \\
   \And
 Jingsai Liang  \\
Department of Computer Science and Data Science\\
Westminster University,\\
Salt Lake City, UT 84105\\
  \texttt{liang@westminsteru.edu} \\
}
\begin{document}
\maketitle
\begin{abstract}
\indent This paper presents an improved physics-informed neural network for simulating the spatio-temporal solution profile of the multidimensional coupled Burgers' equations with high Reynolds numbers. As time evolves, the sharp shock fronts emerge in the solution, creating significant computational challenges for the conventional mesh-based numerical methods. In particular, numerical methods such as finite differences and finite elements suffer from poor stability and strong mesh dependency when resolving the steep solution gradients. To address these challenges, the proposed framework employs a layerwise self-adaptive weighting strategy that dynamically adjusts the penalty weights for the physics residual, initial conditions, and boundary conditions throughout training. Moreover, the framework uses a dual-phase optimization strategy to achieve more stable convergence. 

To check the effectiveness and accuracy of the proposed framework, a set of numerical experiments is conducted to compare it with standard Physics-Informed Neural Network (PINN) with and without Limited-memory Broyden–Fletcher–Goldfarb–Shanno (L-BFGS) optimization. Numerical results exhibit that the proposed framework achieves higher accuracy in terms of relative $L_2-$ error norm than the standard PINN and is able to capture the development of sharp shock fronts as time evolves in the solution.
\end{abstract}


\section{Introduction}
In recent years, research on Neural Networks (NNs) has received a lot of attention from the scientific community due to their huge success in high-impact fields such as generative drug discovery \cite{sb}, speech recognition \cite{hg}, computer vision \cite{ka}, climate forecasting \cite{pj}, natural language translation \cite{wy}, and precision oncology \cite{jj}. This momentum has naturally shifted toward scientific computing, where NNs are used to simulate the spatio-temporal solution profiles of partial differential equations (PDEs) \cite{rm}. However, traditional black-box AI models, such as standard Deep Neural Networks (DNNs), excel at finding complex patterns in massive datasets, but operate without any inherent understanding of the physical laws governing that data \cite{add}. As a consequence, the flexibility of these networks as universal function approximators comes at the cost of a massive parameter space, requiring large volumes of training data that are often unavailable in scientific contexts. Therefore, they often produce physically inconsistent results, and their performance drops sharply when data is sparse or noisy.

To address these drawbacks, Physics-Informed Neural Networks (PINNs) have emerged as a revolutionary framework by embedding the governing PDEs directly into the neural network's loss function \cite{rm}. This approach transforms the network into a ``physics-aware" solver that penalizes any prediction contradicting fundamental principles, enabling accurate simulation even in data-starved environments. Despite this, standard PINNs and conventional mesh-based numerical methods face significant hurdles in simulating the solution profiles of multidimensional coupled Burgers' equations (MCBEs) at high Reynolds numbers. Conventional methods, such as finite element and finite difference methods, rely on extremely fine computational grids, which become computationally prohibitive due to the curse of dimensionality \cite{qa}, and they require specialized shock-capturing techniques, such as Total Variation Diminishing (TVD) or Essentially Non-Oscillatory (ENO) schemes, to prevent oscillations near discontinuities. On the other hand, ``vanilla" PINNs often struggle to solve MCBEs with high Reynolds numbers due to the dominant convective terms creating a highly ``stiff" loss landscape, causing optimizers to stall in local minima and fail to capture the sharp gradients characteristic of the convective terms.

In recent years, standard and improved PINNs have been successfully applied to the 1D and 2D Burgers' equations. To name a few, in \cite{rm}, the authors introduced a standard PINN to solve forward and inverse problems involving the 1D Burgers equation as a benchmark. Jagtap et al. in \cite{ja} proposed a conservative PINN to solve the 1D and 2D Burgers equations. In \cite{mx}, Meng et al. introduced a parareal PINN to solve time-dependent PDEs and applied it to the 1D Burgers equation to demonstrate its performance. In \cite{mrz}, authors developed a sequential method to train PINNs for PDEs sequentially over successive time segments using a single neural network, which solved Cahn-Hilliard equations, Allen Cahn equations, and the 1D parameteric Burgers equation. McClenny and Braga-Neto introduced a novel self-adaptive PINN model that utilizes fully trainable pointwise weights and a soft attention mechanism to solve various PDEs, including 2D coupled Burgers' equations (CBEs) with moderate Reynolds number, in \cite{mcl}. In \cite{orm}, the authors developed a hybrid computational framework that combines PINN and boundary layer theory to solve the 1D Burgers' equations.  Recently, Song et al. in \cite{sy} introduced a loss-attentional PINNS to improve PINN accuracy and convergence by utilizing a Loss-Attentional Network (LAN) to dynamically weight training points and successfully implemented it to 2D CBEs with moderate Reynolds number. In \cite{xw}, the manuscript introduced a weak-form PINN to solve the 1D Burgers equation, by utilizing an integral formulation of the PDE loss to improve training stability in regions with shocks or steep gradients, where standard PINNs often fail. Zhang et al. in \cite{zcz} introduced two improved PINNs, one is a sequential PINN and the other is PINNs with the WENO indicator to solve 1D Burgers' equations with vanishingly small viscosity. 

Most of the existing studies have successfully applied various PINNs models to the 1D Burgers' equation at moderate Reynolds numbers, with implementations primarily in Python. To the best of the author's knowledge, no studies have been conducted to solve 2D and 3D CBEs with high Reynolds numbers using the PINNs models. Therefore, the main objective of this study is to advance the field by introducing an improved PINNs framework incorporating a layerwise self-adaptive weighting strategy and dual-optimization strategy to accurately predict the spatio-temporal dynamics of the following $d$-dimensional CBEs for the velocity vector $\mathbf{u}(\mathbf{x}, t) = [u_1, u_2, \dots, u_d]^T$ with high Reynolds numbers given by: 
\begin{equation}
\begin{cases} 
\displaystyle \frac{\partial u_i}{\partial t} + \sum_{j=1}^d u_j \frac{\partial u_i}{\partial x_j} = \nu \sum_{j=1}^d \frac{\partial^2 u_i}{\partial x_j^2}, & \mathbf{x} \in \Omega, \ t \in (0, T], \ i = 1, \dots, d \\ \\
u_i(\mathbf{x}, 0) = u_{i,0}(\mathbf{x}), & \mathbf{x} \in \Omega, \ i = 1, \dots, d \\ \\
Bu_i(\mathbf{x}, t) = g_i(\mathbf{x}, t), & \mathbf{x} \in \partial\Omega, \ t \in (0, T], \ i = 1, \dots, d
\end{cases}
\label{eq1}
\end{equation}
where \(\nu\) denotes the viscosity parameter, which corresponds to the inverse of the Reynolds number Re, that is, \( \nu = \frac{1}{\text{Re}}\), $d=1,2,3$, $B$ is the type of boundary condition considered (usually Dirichlet, Neumann, or mixed), $u_{i,0}$, and $g_i$ are initial and boundary data, respectively. 

The effectiveness and accuracy of the proposed framework will be validated by comparing it with classic PINNs, with and without Limited-memory Broyden–Fletcher–Goldfarb–Shanno (L-BFGS) optimization using relative $L_2-$error norms. Moreover, the framework will be implemented in MATLAB using its DL toolbox for automatic differentiation and neural network construction.

The following nomenclature is adopted for the frameworks considered in this study to carry out the comparison among the different optimization strategies:
\begin{itemize}
    \item LSAAL-PINN: The layerwise self-adaptive PINN trained with Adam and L-BFGS optimizers.
    \item AL-PINN: The standard PINN trained with Adam and L-BFGS optimizers.
    \item A-PINN: The standard PINN trained only with Adam optimizer.
\end{itemize}

The remainder of the manuscript is organized as follows. In Section 2, standard PINN, along with a layerwise self-adaptive PINN model, are described in detail as the solution approximator of the 3D Burgers' equation. In Section 3, numerical experiments are conducted to show the performance advantage of the proposed PINN framework over the standard PINN framework. Finally, the conclusions are presented in Section 4.

\section{Layerwise Self-Adaptive PINN Trained with Adam and L-BFGS Optimizers (LSAAL-PINN)}
This section provides a detailed overview of the proposed LSAAL-PINN to solve the 3D scalar Burgers' equation given by:
\begin{equation}
\begin{cases} 
\displaystyle u_t + u(u_x+u_y+u_z)=\nu(u_{xx}+u_{yy}+u_{zz}), & \mathbf{x}=(x,y,z) \in \Omega\subset \mathbb{R}^{d=3}, \ t \in (0, T] \\ \\
u(\mathbf{x}, 0) = u_{i,0}(\mathbf{x}), & \mathbf{x} \in \Omega\\ \\
u(\mathbf{x}, t) = g(\mathbf{x}, t), & \mathbf{x} \in \partial\Omega, \ t \in (0, T].
\end{cases}
\label{eq2}
\end{equation}

To solve the 3D Burgers' equation, the solution $u(\mathbf{x},t)$ is approximated by DNN, defined as a continuous and differentiable composite mapping $\Phi_\theta:\mathbb{R}^{d+1}\rightarrow \mathbb{R}$ parametrized by trainable parameters $\theta$ such that $u(\mathbf{x},t)\approx\hat{u}=\Phi_\theta(\mathbf{x}, t)$.

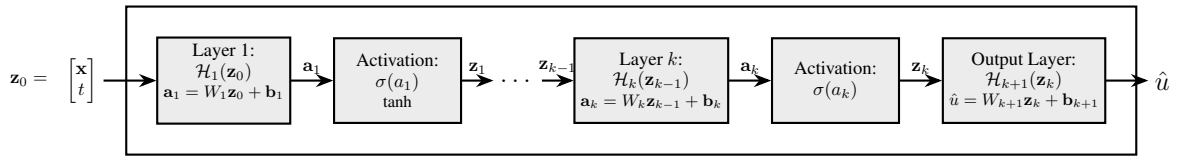
\begin{figure}[H]
    \centering
    \begin{tikzpicture}[
        scale=0.7, transform shape,
        node distance=0.6cm and 0.8cm,
        block/.style={
            rectangle, 
            draw, 
            fill=gray!15, 
            minimum width=2.4cm, 
            minimum height=1.6cm, 
            align=center, 
            thick
        },
        container/.style={
            draw, 
            rectangle, 
            inner sep=0.4cm, 
            thick
        },
        arrow/.style={-Stealth, thick}
    ]

    \node (z0_label) at (-1, 0) {$\mathbf{z}_0 = $};
    \node (coords) [right=0.1cm of z0_label] {
        $\begin{bmatrix} \mathbf{x} \\ t \end{bmatrix}$
    };
 
    \node[block, right=1cm of coords] (L1) {
        Layer 1:\\ 
        $\mathcal{H}_1(\mathbf{z}_0)$\\ 
        \small $\mathbf{a}_1 = W_1 \mathbf{z}_0 + \mathbf{b}_1$\\
    };

    \node[block, right=of L1] (A1) {
        Activation:\\ 
        $\sigma(a_1)$\\ 
        \small tanh
    };

    \node[right=0.6cm of A1] (dots) {\Large $\dots$};

    \node[block, right=0.6cm of dots] (Lk) {
        Layer $k$:\\ 
        $\mathcal{H}_k(\mathbf{z}_{k-1})$\\ 
        \small $\mathbf{a}_k = W_k \mathbf{z}_{k-1} + \mathbf{b}_k$
    };
    \node[block, right=of Lk] (Ak) {
        Activation:\\ 
        $\sigma(a_k)$
    };

    \node[block, right=of Ak] (Out) {
        Output Layer:\\ 
        $\mathcal{H}_{k+1}(\mathbf{z}_k)$\\ 
        \small $\hat{u} = W_{k+1} \mathbf{z}_k + \mathbf{b}_{k+1}$
    };

    \node[right=0.8cm of Out] (u) {\Large $\hat{u}$};
    \draw[arrow] (coords) -- (L1);
    \draw[arrow] (L1) -- node[above] {$\mathbf{a}_1$} (A1);
    \draw[arrow] (A1) -- node[above] {$\mathbf{z}_1$} (dots);
    \draw[arrow] (dots) -- node[above] {$\mathbf{z}_{k-1}$} (Lk);
    \draw[arrow] (Lk) -- node[above] {$\mathbf{a}_k$} (Ak);
    \draw[arrow] (Ak) -- node[above] {$\mathbf{z}_k$} (Out);
    \draw[arrow] (Out) -- (u);
    \begin{scope}[on background layer]
        \node[container, fit=(L1) (Out) ] (box) {};
    \end{scope}

    \end{tikzpicture}
    \caption{DNN $\Phi_\theta :\mathbb{R}^{d+1}\rightarrow \mathbb{R}$}
    \label{fig:figNN}
\end{figure}

For an input $\mathbf{z_0}=(\mathbf{x}, t)^\mathsf{T}$, the feedforward NN with $k$ hidden layers and $n$ neurons in each, defined as a sequence of nested affine transformations and nonlinear activation functions as shown in Figure \ref{fig:figNN}: 
\[\hat{u}=\Phi_\theta(\mathbf{z_0})=\big(\mathcal{H}_{k+1} \circ \sigma \circ \mathcal{H}_k \circ \cdots \circ \sigma \circ \mathcal{H}_1 \big)(\mathbf{z_0}),\]

where 
\begin{itemize}
    \item First Hidden Layer 1: The initial mapping $\textbf{a}_1=\mathcal{H}_1(\mathbf{z}_{0}) = \mathbf{W}_1 \mathbf{z}_{0}+\mathbf{b}_1,~\text{where}~ \mathbf{W}_1 \in \mathbb{R}^{n\times (d+1)}, \mathbf{b}_1\in \mathbb{R}^n$ project the input $\mathbf{z}_0$ into a high-dimensional feature space of width $n$.
    \item Intermediate Hidden Layers: The network processes the features through $k-1$ additional hidden layers  $\mathbf{a}_i=\mathcal{H}_i(\mathbf{z}_{i-1}) = \mathbf{W}_i \mathbf{z}_{i-1}+\mathbf{b}_i, \mathbf{W}_i\in \mathbb{R}^{n\times n}, \mathbf{b}_i\in \mathbb{R}^n,~i = 2, 3, \cdots,k$. These layers enable the network to understand complex nonlinear phenomena.
     \item Activation (nonlinear processing): For each hidden layer, a nonlinear activation function $\textbf{z}_i=\sigma(\textbf{a}_i)$ is applied. The activation function $\sigma=\tanh$ is then sandwiched between affine layers to provide nonlinearity and smoothness required to approximate PDE solutions and compute derivatives required in the physics residual.
    \item Output Layer: $\hat{u}=\mathcal{H}_{k+1}(\mathbf{z}_k)=\mathbf{W}_{k+1}\mathbf{z}_k+\mathbf{b}_{k+1}, \mathbf{W}_{k+1}\in \mathbb{R}^{1\times n}, \mathbf{b}_{k+1}\in \mathbb{R}.$ 
    \item $\theta = \left\{ \mathbf{W}_{i},\mathbf{b}_i\right\}_{i=1}^{k+1}$ is the set of trainable weights and biases. Throughout this study, trainable parameters $\theta$ are initialized using the Xavier initialization scheme \cite{xyb} to prevent exploding and vanishing gradients in $\tanh$ networks.

\end{itemize}
To ensure that $\hat{u}$ gradually converges to the true solution $u(\mathbf{x}, t)$, $\hat{u}$ is trained to satisfy the physical model by embedding the equations in Eq. \eqref{eq2} into the optimization process through the total loss function. To embed the underlying physics into the NN, the Automatic Differentiation (AD) is utilized to compute the partial derivatives: \[ \hat{u}_t, \hat{u}_x, \hat{u}_y, \hat{u}_z, \hat{u}_{xx}, \hat{u}_{yy},~\text{and}~ \hat{u}_{zz}\] by repeated chain rule through the composite neural map which yields the physics residual function as:
\[f:=\hat{u}_t + \hat{u}(\hat{u}_x+\hat{u}_y+\hat{u}_z)-\nu(\hat{u}_{xx}+\hat{u}_{yy}+\hat{u}_{zz}).\] The exact solution $u(\mathbf{x}, t)$ satisfies $f=0$ for all $(x,y,z,t)$ and therefore the NN is trained to minimize the total loss function over the given domain. The total loss function simultaneously penalizes the following three types of residuals.

Let $\left\{\textbf{x}_i, t_i\right\}_{i=1}^{N_{col}}$ represent a set of collocation points sampled in the interior of the domain. The PDE residual is then defined as: 
\[ \mathcal{L}_{\text{PDE}}=\frac{1}{N_{col}}\sum_{i=1}^{N_{col}}|f|^2,\]

similarly, let $\left\{\textbf{x}_j, 0, u_{0_{j}}\right\}_{j=1}^{N_{ic}}$ be the initial points, and $\left\{\textbf{x}_k, t_k,g_k\right\}_{k=1}^{N_{bc}}$ be the boundary points. The initial and boundary residuals are then defined as:
\[ \mathcal{L}_{\text{IC}}=\frac{1}{N_{ic}}\sum_{j=1}^{N_{ic}}|\hat{u}(\textbf{x}_j,0)-u_{0_{j}}|^2,\] 

\[ \mathcal{L}_{\text{BC}}=\frac{1}{N_{bc}}\sum_{k=1}^{N_{ic}}|\hat{u}(\textbf{x}_k,t_k)-g(\textbf{x}_k,t_k)|^2.\]

Combining the above three residuals, the total loss function is formulated as a minmax $\displaystyle(\min_{\theta}\max_{\lambda}\mathcal{L}(\theta, \lambda))$ optimization problem:
\[\mathcal{L}(\theta, \lambda)
=
\lambda_{\text{PDE}}\mathcal{L}_{\text{PDE}}(\theta)
+
\lambda_{IC}\mathcal{L}_{IC}(\theta)
+
\lambda_{BC}\mathcal{L}_{BC}(\theta),\]
where $\lambda=[\lambda_{\text{PDE},~ \lambda_{IC}, ~\lambda_{BC}}]$ are self adaptive weights related to each constraint and $\theta$ are the NN parameters. In the work, the training procedure uses a dual-phase optimization approach to achieve a stable solution strategy for the minmax learning problem. In phase one, the parameters $\theta$ and $\lambda$ are updated iteratively using a gradient-based scheme over a fixed number of iterations. At each iteration, the loss function $\mathcal{L}(\theta, \lambda)$, along with its gradients $\nabla_{\theta}\mathcal{L},~\nabla_{\lambda}\mathcal{L}$ with respect to the $\theta$ and $\lambda$, respectively are evaluated using AD. Once the gradients are obtained, the parameters $\theta$ are updated using the Adam optimizer, leading to updates of the form:\[\theta^{k+1}=\text{Adam}(\theta^k, \nabla_{\theta}\mathcal{L} ).\]

Simultaneously, adaptive weights $\lambda$ are updated via the Gradient Ascent strategy to maximize the total loss with respect to these variables. The update rule is derived by reversing the direction of the gradient of the loss function with respect to $\lambda$, resulting in updates of the form: 
\[\lambda^{k+1}=\lambda^k+\eta_\lambda\nabla_{\lambda}\mathcal{L}(\theta, \lambda),\] 
where $\eta_\lambda$ is the learning rate, to dynamically determine the relative importance of each loss component during optimization. This self-adaptive weighting strategy ensures that no single residual term dominates training, if it happens to be larger than the others. To maintain numerical stability, these parameters are mapped through a logarithmic transformation $$\lambda_j=\log (\lambda_j-\min(\lambda)+1)+\epsilon,~~~j\in \left\{\text{PDE, IC, BC}\right\}$$ with a small constant $\epsilon>0$ to ensure strict positivity.

After completion of the Adam-based training phase, in phase two, the adaptive weights $\lambda$ are frozen, and the optimization transitions to a second refinement stage using the second-order optimizer L-BFGS algorithm to accelerate the convergence. This optimizer utilizes the optimized parameters $\theta$ obtained via the first-order optimization phase and refines the DNN mapping $\Phi_{\theta}$ by solving the reduced problem:\[\min_{\theta}\mathcal{L}(\theta, \lambda)\]
until it satisfies the Burger's equation within a prescribed tolerance. As a result, the proposed PINN framework yields a mesh-free and continuous representation of the solution profile that satisfies the governing equation together with the initial and boundary constraints. 

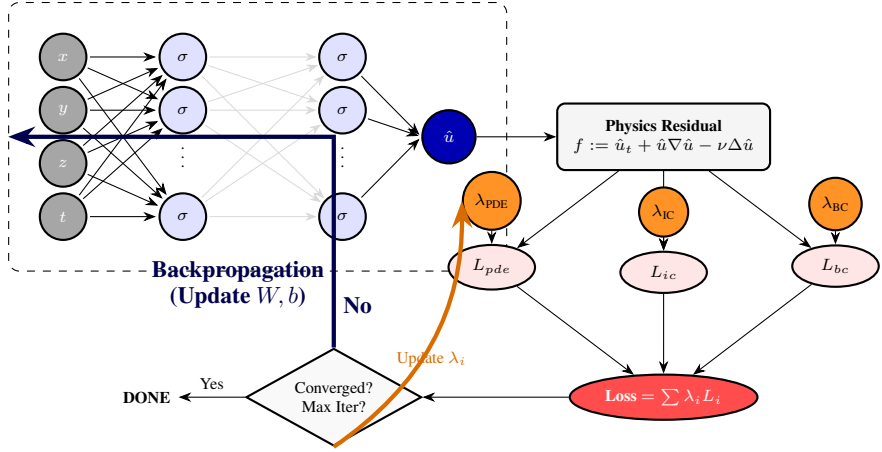
\begin{figure}[H]
\begin{tikzpicture}[
    >=Stealth,
    scale=0.88,
    transform shape,
    neuron/.style={circle, draw, thick, fill=blue!12, minimum size=7mm},
    input/.style={circle, draw, thick, fill=gray!70, text=white, minimum size=7mm},
    output/.style={circle, draw, thick, fill=blue!70!black, text=white, minimum size=8mm},
    loss/.style={ellipse, draw, thick, fill=red!10, minimum width=1.3cm},
    weight/.style={circle, draw, thick, fill=orange!85, minimum size=6mm},
    block/.style={rectangle, draw, thick, rounded corners=3pt, fill=gray!8, minimum width=3.2cm, minimum height=1cm, align=center},
    total/.style={ellipse, draw, thick, fill=red!70, text=white, minimum width=2.4cm},
    decision/.style={diamond, draw, thick, aspect=1.8, fill=gray!5, align=center},
    every node/.style={font=\scriptsize}
]

\node[input] (x) at (0,1.2) {$x$};
\node[input] (y) at (0,0.4) {$y$};
\node[input] (z) at (0,-0.4) {$z$};
\node[input] (t) at (0,-1.2) {$t$};

\node[neuron] (h1-1) at (1.8,1.2) {$\sigma$};
\node[neuron] (h1-2) at (1.8,0.4) {$\sigma$};
\node (h1-dots) at (1.8,-0.2) {$\vdots$};
\node[neuron] (h1-n) at (1.8,-1.2) {$\sigma$};

\node (layer-dots) at (3,0) {$\cdots$};

\node[neuron] (hk-1) at (4.2,1.2) {$\sigma$};
\node[neuron] (hk-2) at (4.2,0.4) {$\sigma$};
\node (hk-dots) at (4.2,-0.2) {$\vdots$};
\node[neuron] (hk-n) at (4.2,-1.2) {$\sigma$};

\node[output] (u) at (5.8,0) {$\hat{u}$};

\foreach \i in {x,y,z,t}
    \foreach \j in {1,2,n}
        \draw[->, thin, shorten >=1pt, shorten <=1pt] (\i) -- (h1-\j);

\foreach \i in {1,2,n}
    \foreach \j in {1,2,n}
        \draw[->, ultra thin, black!15, shorten >=1pt, shorten <=1pt] (h1-\i) -- (hk-\j);

\foreach \i in {1,2,n}
    \draw[->, thin, shorten >=1pt, shorten <=1pt] (hk-\i) -- (u.west);

\node[draw,dashed,rounded corners,inner sep=4mm,
fit=(x)(t)(hk-1)(u),
] (fitbox) {};

\node[block, right=1.2cm of u] (res) {
\textbf{Physics Residual}\\
$f:=\hat u_t+\hat u\nabla \hat u-\nu\Delta \hat u$
};
\draw[->] (u)--(res);

\node[loss, below left=1.2cm and 0.5cm of res] (Lpde) {$L_{pde}$};
\node[loss, below=1.2cm of res] (Lic) {$L_{ic}$};
\node[loss, below right=1.2cm and 0.5cm of res] (Lbc) {$L_{bc}$};

\draw[->] (res)--(Lpde);
\draw[->] (res)--(Lic);
\draw[->] (res)--(Lbc);

\node[weight, above=2mm of Lpde] (wpde) {$\lambda_{\text{PDE}}$};
\node[weight, above=2mm of Lic] (wic) {$\lambda_{\text{IC}}$};
\node[weight, above=2mm of Lbc] (wbc) {$\lambda_{\text{BC}}$};

\draw[->, dashed] (wpde)--(Lpde);
\draw[->, dashed] (wic)--(Lic);
\draw[->, dashed] (wbc)--(Lbc);

\node[total, below=1.2cm of Lic] (total)
{\textbf{Loss} $= \sum \lambda_i L_i$};

\draw[->] (Lpde)--(total);
\draw[->] (Lic)--(total);
\draw[->] (Lbc)--(total);

\node[decision, left=2.2cm of total] (check) {Converged?\\Max Iter?};
\draw[->] (total.west) -- (check.east);

\node[left=1cm of check] (done) {\scriptsize\textbf{DONE}};
\draw[->] (check.west) -- node[above]{Yes} (done.east);

\draw[->, ultra thick, black!70!blue]
(check.north) |- 
node[pos=0.1, right, font=\bfseries] {No}
node[pos=0.15, left, align=center, font=\bfseries] {Backpropagation\\(Update $W,b$)}
(fitbox.west);

\draw[->, ultra thick, orange!85!black, bend right=28]
(check.south) to node[midway, below]{\scriptsize Update $\lambda_i$} (wpde.west);

\end{tikzpicture}
\caption{The schematics of the proposed LSAAL-PINN for the 3D Burgers' equation}
\end{figure}

\section{Numerical Experiments and Discussion}
In this section, three test problems with high Reynolds numbers are considered to compare the effectiveness and accuracy of the proposed framework LSAAL-PINN. To conduct a comparative analysis, the LSAAL-PINN is compared against the AL-PINN and A-PINN. The accuracy of the methods is measured in terms of relative $L_2-$ error norm defined as $\frac{\|u-\hat{u}\|_2}{\|u\|_2}$, where $\|\cdot\|_2$ represents $L_2$ norm.

\begin{figure}[H]
\begin{minipage}[b]{0.33\linewidth} 
\centering
\centerline{\includegraphics[width=\linewidth,height=\textheight,keepaspectratio]{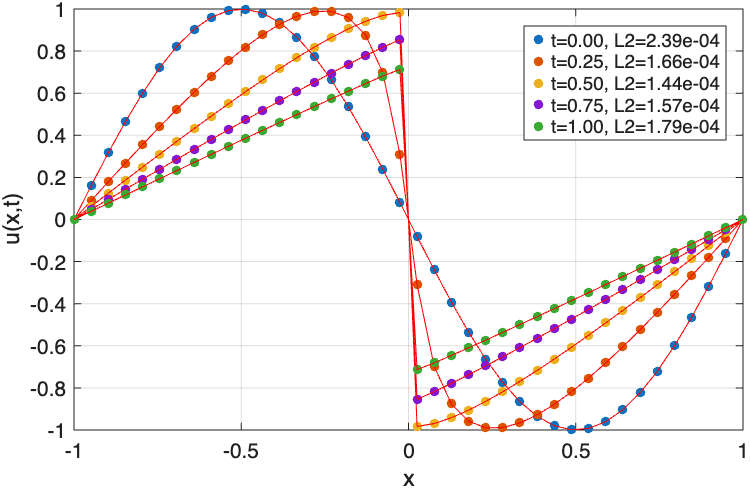}}
\centerline{\text{(a) LSAAL-PINN}}
\end{minipage}
\begin{minipage}[b]{0.33\linewidth} 
\centering
\centerline{\includegraphics[width=\linewidth,height=\textheight,keepaspectratio]{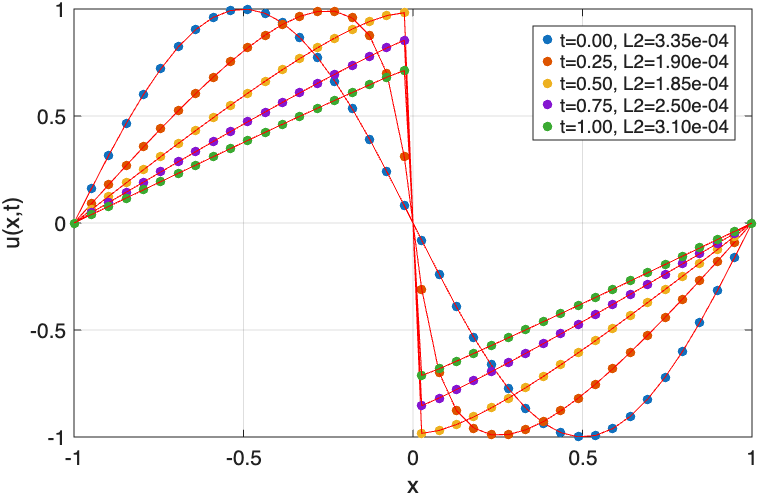}}
\centerline{\text{(b) AL-PINN}}
\end{minipage}
\begin{minipage}[b]{0.33\linewidth} 
\centering
\centerline{\includegraphics[width=\linewidth,height=\textheight,keepaspectratio]{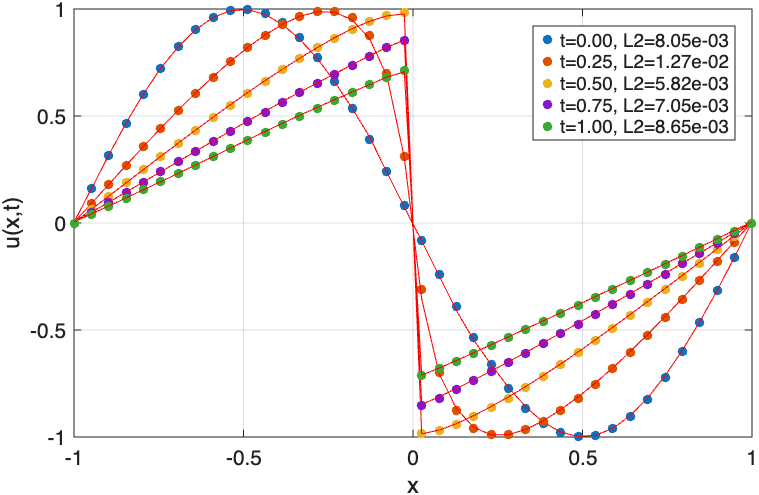}}
\centerline{\text{(c) A-PINN}}
\end{minipage}
\caption{\footnotesize{The comparison of exact solution (dots) vs. predicted solution (solid lines) $u(x,t)$ at various time levels for 1D Burgers' equation. Relative $L_2-$errors for each temporal snapshot are included in the legend. }}
\label{fig: 1D burger}
\end{figure}

\begin{figure}[H]
\begin{minipage}[b]{0.24\linewidth} 
\centering
\centerline{\includegraphics[width=\linewidth,height=\textheight,keepaspectratio]{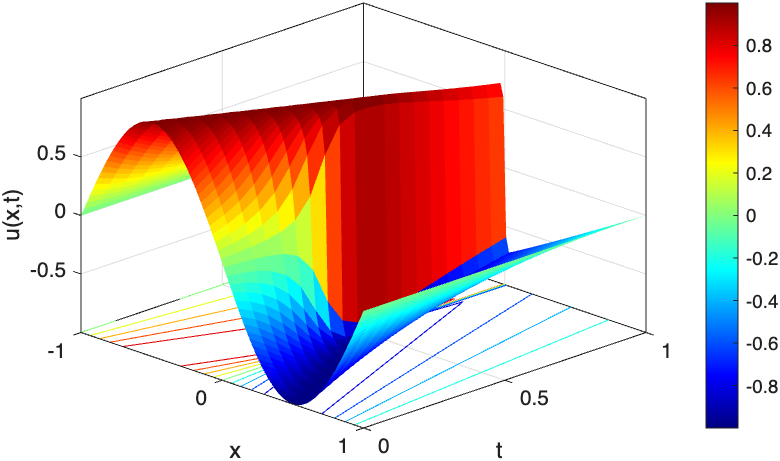}}
\centerline{\text{(a) LSAAL-PINN}}
\end{minipage}
\begin{minipage}[b]{0.24\linewidth} 
\centering
\centerline{\includegraphics[width=\linewidth,height=\textheight,keepaspectratio]{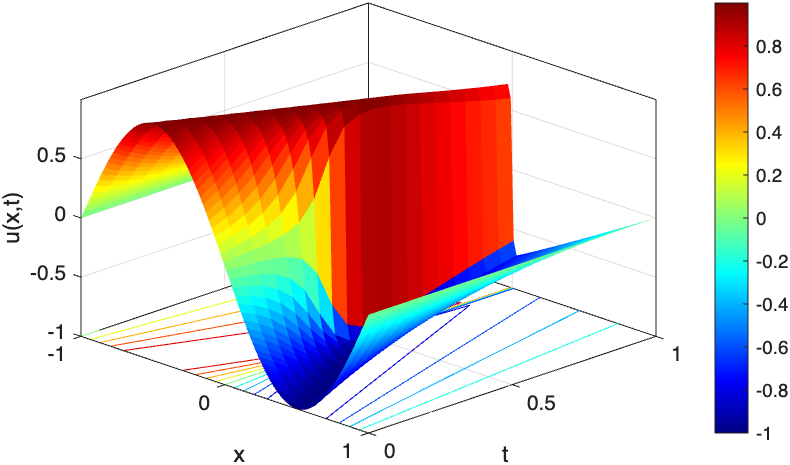}}
\centerline{\text{(b) AL-PINN}}
\end{minipage}
\begin{minipage}[b]{0.24\linewidth} 
\centering
\centerline{\includegraphics[width=\linewidth,height=\textheight,keepaspectratio]{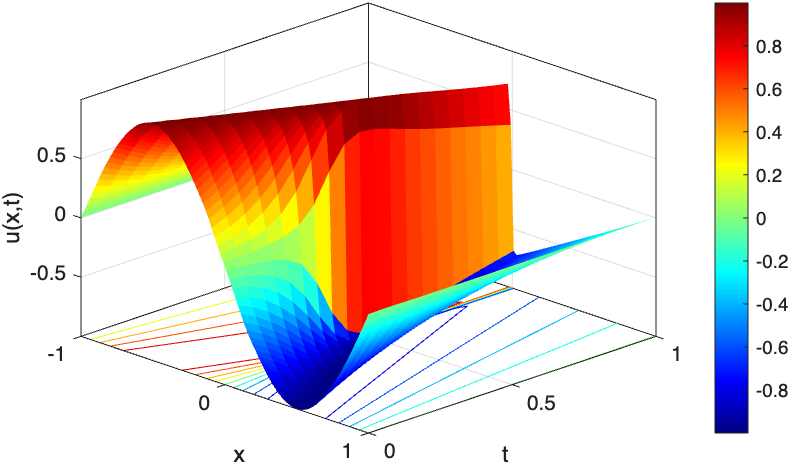}}
\centerline{\text{(c) A-PINN}}
\end{minipage}
\begin{minipage}[b]{0.24\linewidth} 
\centering
\centerline{\includegraphics[width=\linewidth,height=\textheight,keepaspectratio]{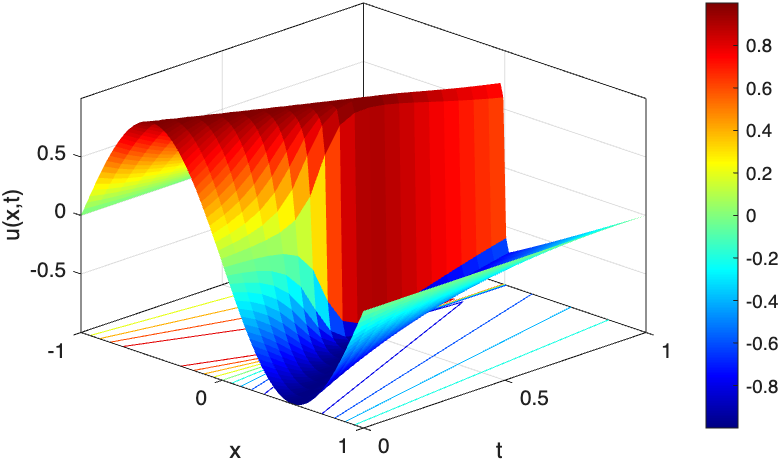}}
\centerline{\text{(d) Exact}}
\end{minipage}
\caption{\footnotesize{The comparison between exact and predicted solution of 1D Burgers' equation}}
\label{fig: 1D Burger surf}
\end{figure}

\begin{figure}[H]
\begin{minipage}[b]{0.33\linewidth} 
\centering
\centerline{\includegraphics[width=\linewidth,height=\textheight,keepaspectratio]{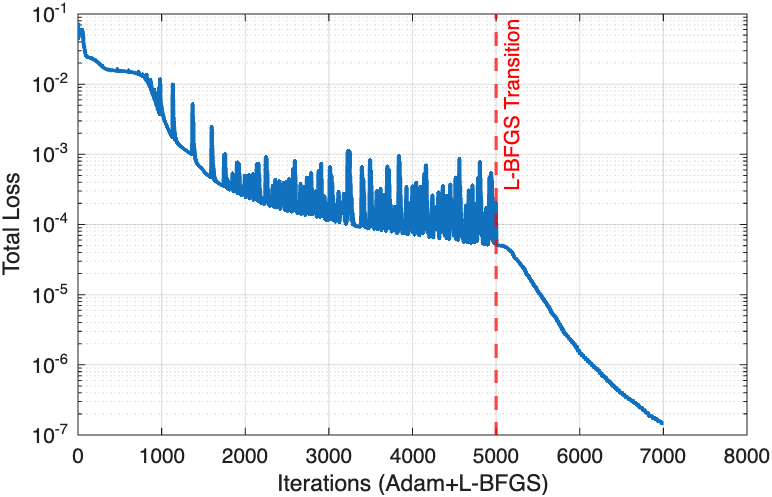}}
\centerline{\text{(a) LSAAL-PINN}}
\end{minipage}
\begin{minipage}[b]{0.33\linewidth} 
\centering
\centerline{\includegraphics[width=\linewidth,height=\textheight,keepaspectratio]{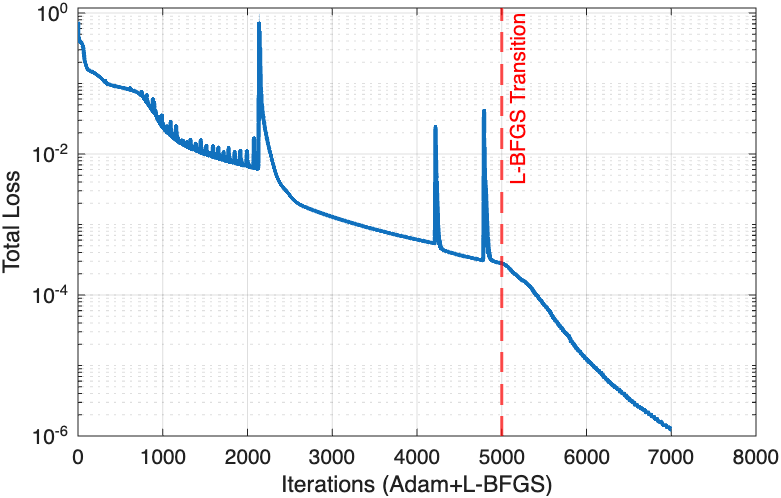}}
\centerline{\text{(b) AL-PINN}}
\end{minipage}
\begin{minipage}[b]{0.33\linewidth} 
\centering
\centerline{\includegraphics[width=\linewidth,height=\textheight,keepaspectratio]{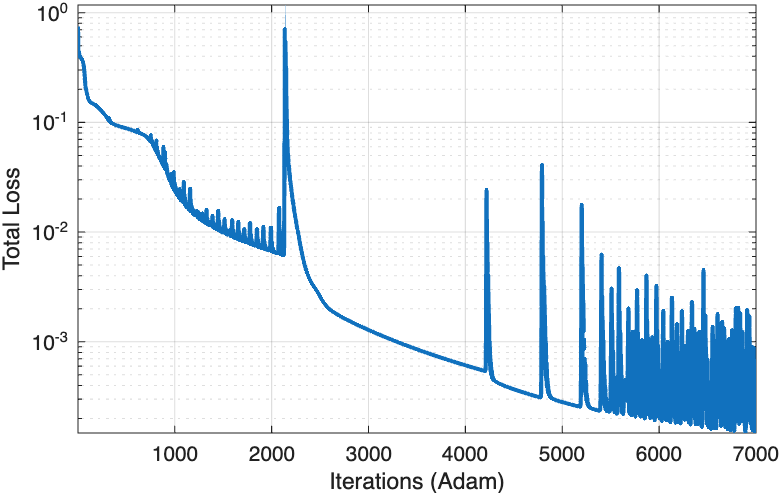}}
\centerline{\text{(c) A-PINN}}
\end{minipage}
\caption{\footnotesize{The comparison of the total training loss values of three PINN variants for Example 1.}}
\label{fig: 1dloss}
\end{figure}

Throughout this study, the initial weights of the DNN are initialized via the Xavier initialization scheme to ensure the stable gradient flow. The hyperbolic tangent $(\tanh)$ is utilized as the activation function across all the hidden layers for the calculation of higher-order physics residuals. To generate well-distributed collocation points across the computational domain, the Latin Hypercube Sampling (LHS) \cite{mst} is utilized over the random sampling. This sampling strategy helps the NN to capture the steep gradient in the solution of Burgers' equation with high Reynold numbers more effectively by ensuring that the points are always present near the transition regions. All training processes utilized a multi-stage optimization strategy, which includes the Adam algorithm for global optimization, followed by the L-BFGS algorithm for high-precision refinement of the loss landscape. Moreover, the following parameters are used to configure each PINN variant to provide clarity and consistency in the performance comparisons presented for the test problems considered below:
\begin{itemize}
    \item Depth and Width: The number of hidden layers and the number of neurons in each hidden layer of the fully connected NN.
    \item LR: The learning rate used in the Adam optimizer.
    \item $N_{\text{col}}$: The size of the data sets for PDE residuals.
    \item $N_{\text{ic}}$: The size of the data sets for initial conditions.
    \item $N_{\text{bc}}$: The size of the data sets for boundary conditions.
    \item Adam Iter: The number of Adam iterations.
    \item L-BFGS Iter: The number of L-BFGS iterations.
\end{itemize}
All code in this work is implemented in MATLAB 2025a using its deep learning toolbox for automatic differentiation. All simulations were performed on a laptop equipped with an Apple M3 Pro chip and 18 GB of unified memory, running macOS Sonoma 14.7.5. 
\begin{example}
Here, we consider the 1D Burgers' equation defined on the domain $x \in [-1, 1]$ and $t > 0$ as a benchmark problem to check the accuracy and effectiveness of the proposed PINN framework:
\begin{flalign*}
    &u_t + u u_x - \nu u_{xx} = 0&,
\end{flalign*}
with the following initial and boundary conditions:
\begin{flalign*}
    &u(x, 0) = -\sin(\pi x)&\\
    &u(-1, t) = u(1, t) = 0,&
\end{flalign*}
where $\nu$ represents the kinematic viscosity coefficient.

The analytical solution of the equation is given by \cite{cmp}:
\begin{flalign*}
&u(x, t) = \frac{-\int_{-\infty}^{+\infty} \sin \pi(x-\eta) f(x-\eta) e^{\frac{-\eta^2}{4\nu t}}  d \eta}{\int_{-\infty}^{+\infty} f(x-\eta) e^{\frac{-\eta^2}{4\nu t}} d \eta},&
\end{flalign*}
where the function $f(y)$ is defined as:
\begin{flalign*}
& f(y) = e^{-\frac{\cos \pi y}{2 \pi \nu}}.&
\end{flalign*}
\end{example}
As detailed in Table \ref{tab:summary}, all three PINN variants utilize a consistent network architecture consisting of 6 hidden layers with 32 neurons in each layer. Moreover, all the frameworks utilize the same hyperparameters and training configurations to approximate the solution profile 
 of the 1D Burgers' equation with $\nu=\frac{0.01}{\pi}$ over the time interval $0\leq t \leq 1.$ As a result of maintaining the consistent network architecture and uniform data sets, the study ensures that any performance variance observed during the experiment is due to the influence of the adopted optimization strategy and proposed self-adaptive weighting strategy. 
\begin{table}[H]
    \centering
    \begin{adjustbox}{width=\textwidth}
    \begin{tabular}{|lcccccccc|}
        \toprule
        \textbf{Method} & \textbf{Depth} & \textbf{Width} & \textbf{LR} & $N_{col}$ & $N_{ic}$ & $N_{bc}$ & \textbf{Adam Iter.} & \textbf{L-BFGS Iter.} \\
        \midrule
        LSAAL-PINN & 6 & 32 & $10^{-3}$ & 5000 & 400 & 200 & 5000 & 2000  \\
        AL-PINN   & 6 & 32 & $10^{-3}$ & 5000 & 400 & 200 & 5000 & 2000 \\
        A-PINN    & 6 & 32 & $10^{-3}$ & 5000 & 400 & 200 & 7000 & -  \\
        \bottomrule
    \end{tabular}
    \end{adjustbox}
    \caption{Summary of hyperparameters and training configurations for the PINN variants used for Example 1.}
    \label{tab:summary}
\end{table}
To compare the accuracy of the three PINN variants, in Figure \ref{fig: 1D burger}, the exact solution vs. predicted solution is plotted until $t=1$. Moreover, relative $L_2-$errors at five different time levels starting from $t=0$ to $t=1$ are computed and depicted in the legend of the Figure \ref{fig: 1D burger}. One can see from Figure \ref{fig: 1D burger} that all three frameworks capture the evolution of a waveform that develops a sharp gradient at $ x=0$. However, the LSAAL-PINN outperforms the rest of the two PINN variants in terms of accuracy by achieving $1.79\mathrm{e}{-4}$ relative $L_2-$ error after $7$k runs, which is better than $3.10\mathrm{e}{-4}$ and $8.65\mathrm{e}{-3}$ obtained from AL-PINN and A-PINN, respectively. Moreover, the Figure \ref{fig: 1D burger} illustrates that the LSAAL-PINN predictions (solid lines) align most precisely with the exact solution (dots), particularly at the shock front. This indicates that the self-adaptive weighting strategy is more effective for capturing the shock profile.  
The spatio-temporal solution profiles of the 1D Burgers' equation predicted via all three PINN variants are compared with the exact solution using the surface plots in Figure \ref{fig: 1D Burger surf} (a)-(d). From the Figure \ref{fig: 1D Burger surf}, we can see that all three variants successfully capture the evolution of the solution and development of the sharp gradient at $x=0$ as expected; however, the degree of accuracy varies significantly. The LSAAL-PINN exhibits the highest precision relative to the exact solution. In contrast, the A-PINN exhibits the most noticeable discrepancies compared to the exact solution because it relies solely on the Adam optimizer. 
In the case of 1D Burgers' equation, the trend of the training convergence behaviors of the three variants is depicted in Figure \ref{fig: 1dloss} throughout the training process. From Figure \ref{fig: 1dloss} (a)-(c), one can see that an initial decay in loss is observed during the first 5000 Adam iterations across all three variants. However, the A-PINN depicts higher loss values and significant oscillations throughout the training period of $7$k runs. In contrast, a sharp decay in total loss is observed after 5000 iterations for the LSAAL-PINN and the AL-PINN, demonstrating the advantage of the hybrid optimization strategy employed by them. Moreover, the total loss of the LSAAL-PINN remains consistently and markedly lower than that of the other two variants throughout the training period, which justifies the importance of the self-adaptive weighting strategy in achieving a high-precision solution of Burgers' equation with high Reynolds number.

\begin{figure}[H]
\begin{minipage}[b]{0.33\linewidth} 
\centering
\centerline{\includegraphics[width=\linewidth,height=\textheight,keepaspectratio]{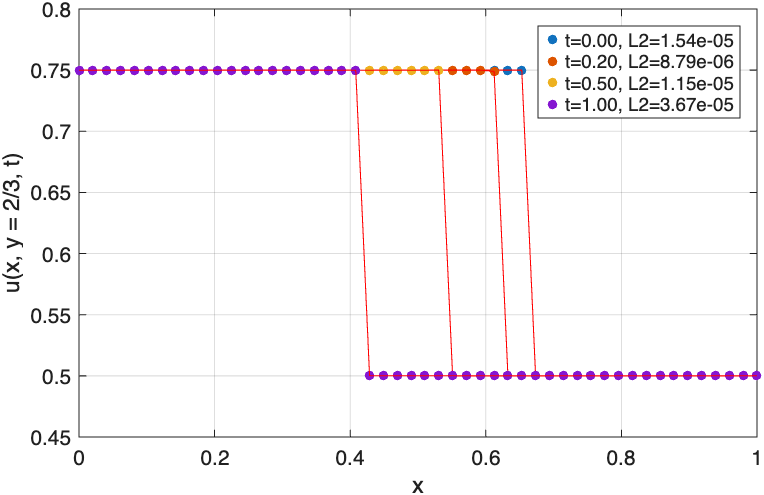}}
\centerline{\text{(a) LSAAL-PINN with $\nu=0.0001$}}
\end{minipage}
\begin{minipage}[b]{0.33\linewidth} 
\centering
\centerline{\includegraphics[width=\linewidth,height=\textheight,keepaspectratio]{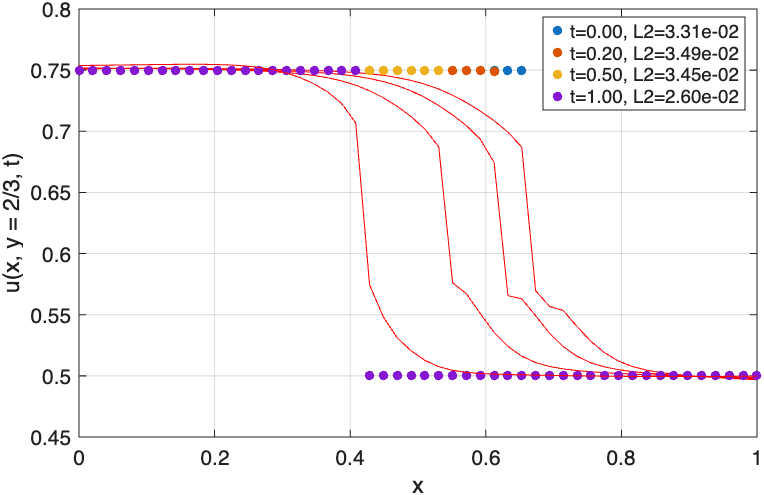}}
\centerline{\text{(b) A-PINN with $\nu=0.0001$}}
\end{minipage}
\begin{minipage}[b]{0.33\linewidth} 
\centering
\centerline{\includegraphics[width=\linewidth,height=\textheight,keepaspectratio]{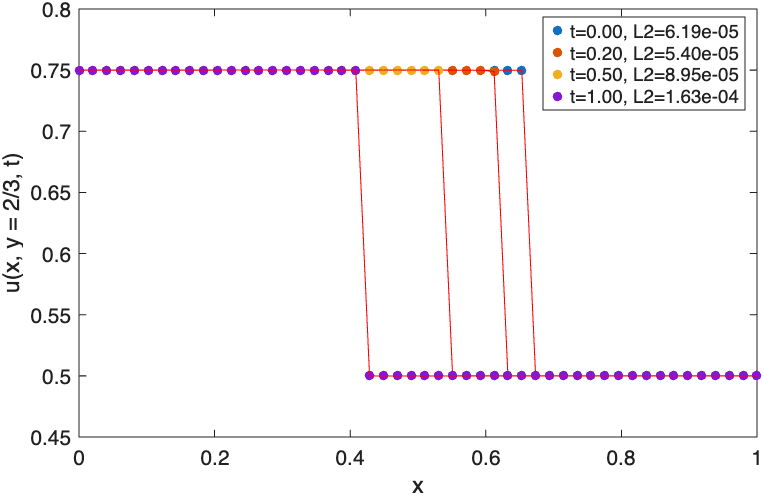}}
\centerline{\text{(c) AL-PINN with $\nu=0.0001$}}
\end{minipage}
\begin{minipage}[b]{0.33\linewidth} 
\centering
\centerline{\includegraphics[width=\linewidth,height=\textheight,keepaspectratio]{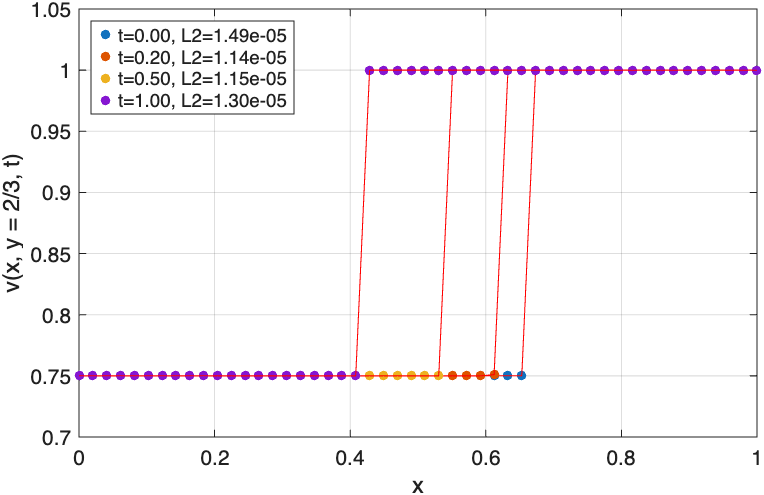}}
\centerline{\text{(d) LSAAL-PINN with $\nu=0.0001$}}
\end{minipage}
\begin{minipage}[b]{0.33\linewidth} 
\centering
\centerline{\includegraphics[width=\linewidth,height=\textheight,keepaspectratio]{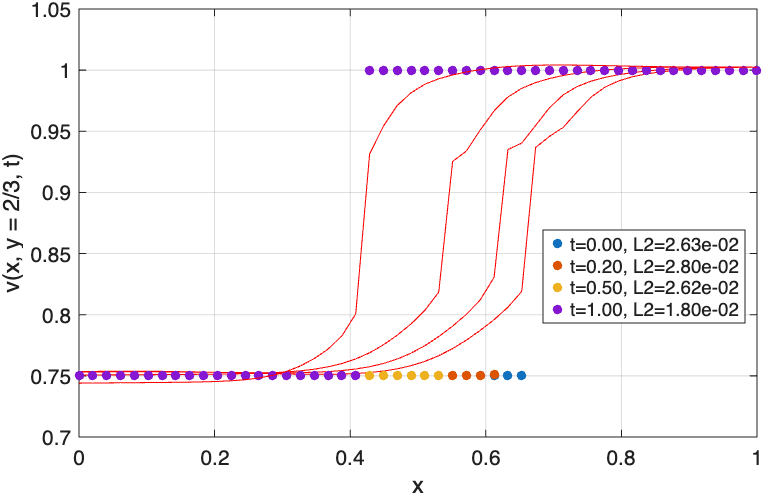}}
\centerline{\text{(e) A-PINN with $\nu=0.0001$}}
\end{minipage}
\begin{minipage}[b]{0.33\linewidth} 
\centering
\centerline{\includegraphics[width=\linewidth,height=\textheight,keepaspectratio]{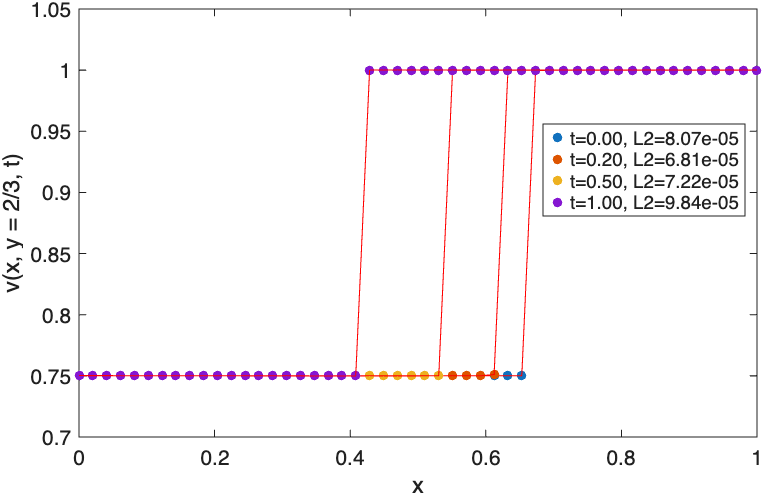}}
\centerline{\text{(f) AL-PINN with $\nu=0.0001$}}
\end{minipage}
\caption{\footnotesize{Comparison of the exact solution (dots) and the predicted solution (solid lines) at various time levels for 2D CBEs. Relative $L_2-$errors for each temporal snapshot are included in the legend. From top to bottom, the first row presents the solution for $u$, while the second row shows the solution for $v$ at $y=\frac{2}{3}$ with high Reynolds number.}}
\label{fig: example2}
\end{figure}

\begin{figure}[H]
\begin{minipage}[b]{0.24\linewidth} 
\centering
\centerline{\includegraphics[width=\linewidth,height=\textheight,keepaspectratio]{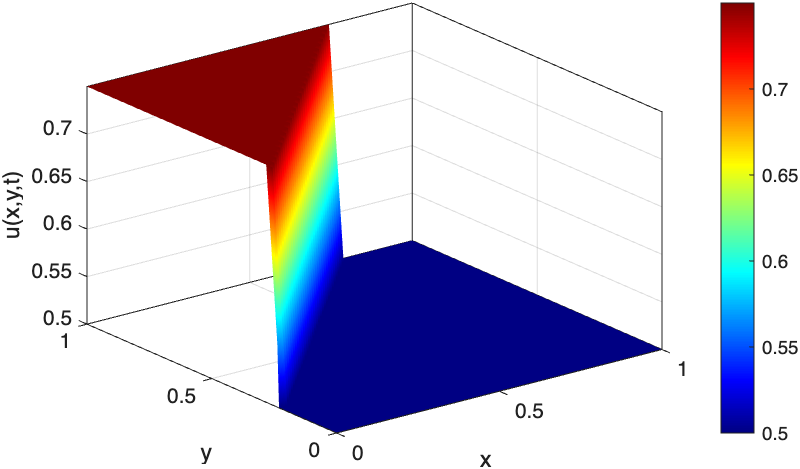}}
\centerline{\text{(a) LSAAL-PINN}}
\end{minipage}
\begin{minipage}[b]{0.24\linewidth} 
\centering
\centerline{\includegraphics[width=\linewidth,height=\textheight,keepaspectratio]{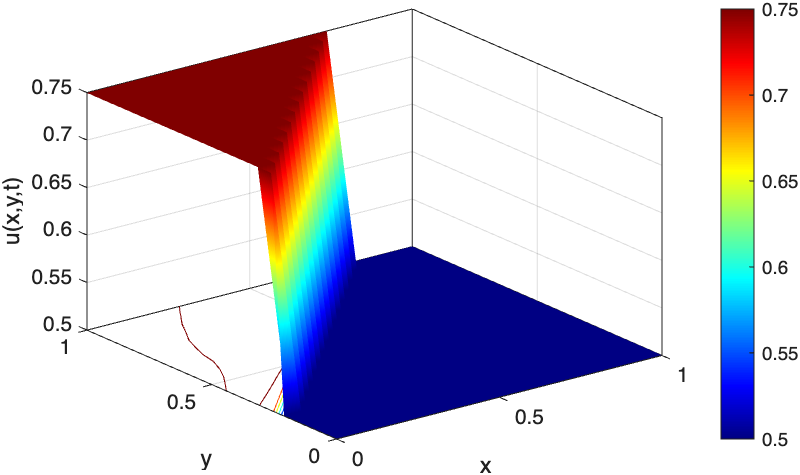}}
\centerline{\text{(b) AL-PINN}}
\end{minipage}
\begin{minipage}[b]{0.24\linewidth} 
\centering
\centerline{\includegraphics[width=\linewidth,height=\textheight,keepaspectratio]{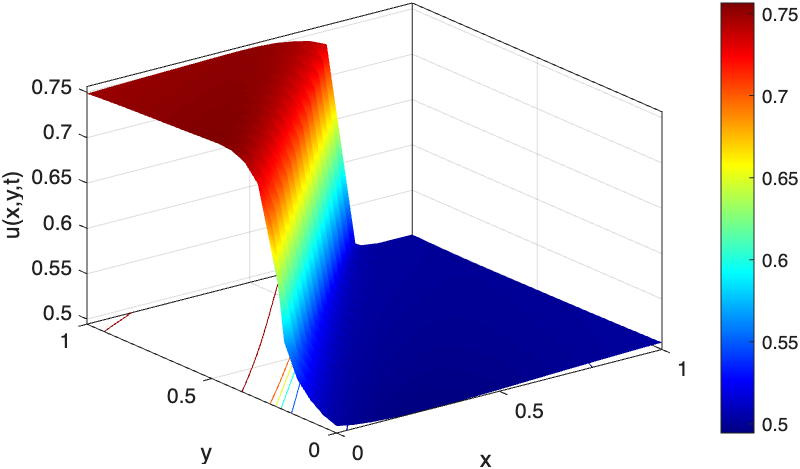}}
\centerline{\text{(c) A-PINN}}
\end{minipage}
\begin{minipage}[b]{0.24\linewidth} 
\centering
\centerline{\includegraphics[width=\linewidth,height=\textheight,keepaspectratio]{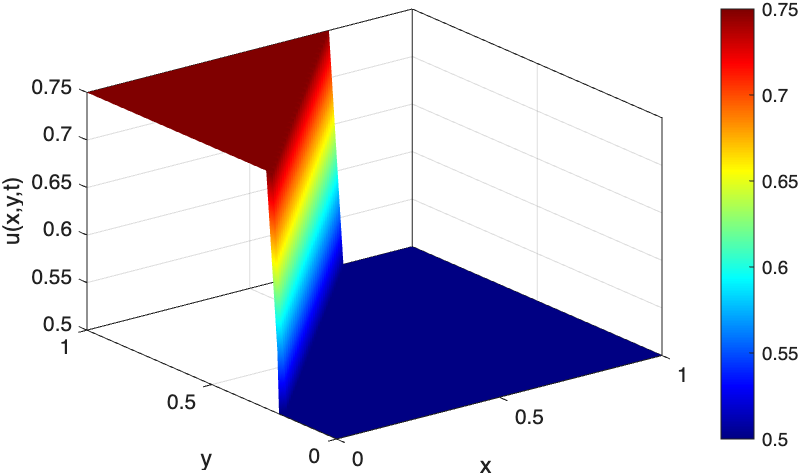}}
\centerline{\text{(d) Exact}}
\end{minipage}
\caption{\footnotesize{The spatio-temporal solution profile of $u$ with $\nu=0.0001$ until time $t=1$ obtained by three PINN variants.}}
\label{fig: 2D Burger surf}
\end{figure}

\begin{figure}[H]
\begin{minipage}[b]{0.33\linewidth} 
\centering
\centerline{\includegraphics[width=\linewidth,height=\textheight,keepaspectratio]{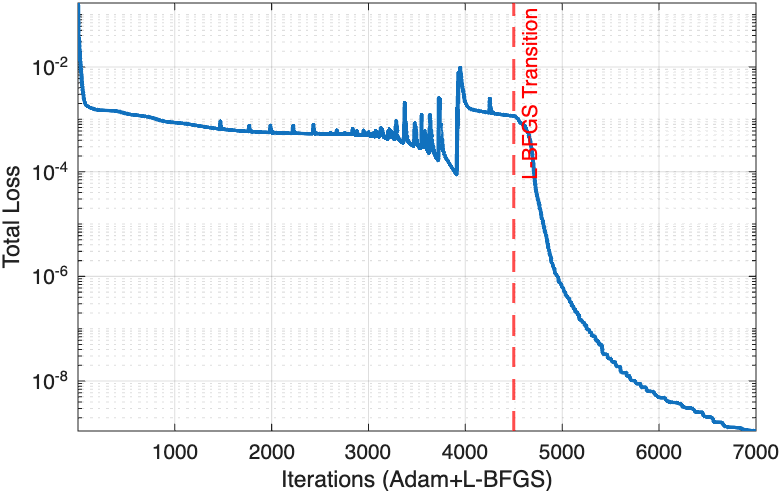}}
\centerline{\text{(a) LSAAL-PINN}}
\end{minipage}
\begin{minipage}[b]{0.33\linewidth} 
\centering
\centerline{\includegraphics[width=\linewidth,height=\textheight,keepaspectratio]{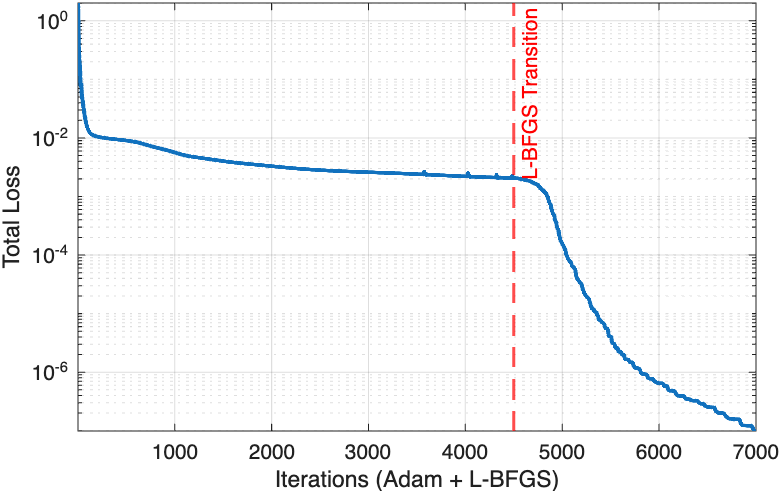}}
\centerline{\text{(b) AL-PINN}}
\end{minipage}
\begin{minipage}[b]{0.33\linewidth} 
\centering
\centerline{\includegraphics[width=\linewidth,height=\textheight,keepaspectratio]{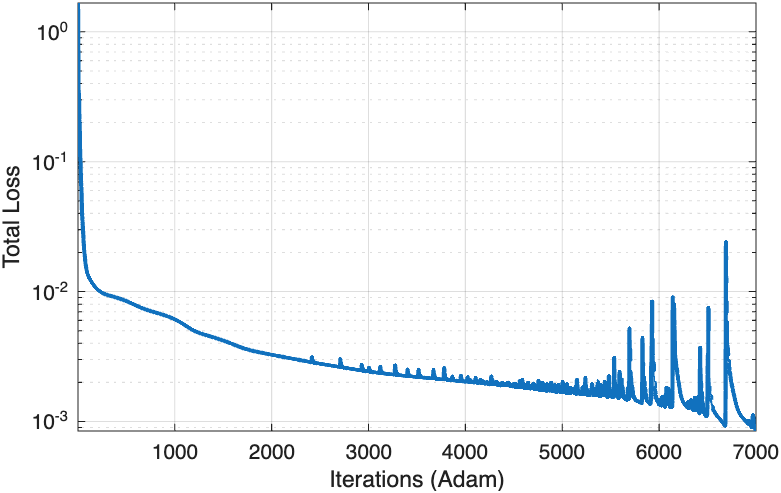}}
\centerline{\text{(c) A-PINN}}
\end{minipage}
\caption{\footnotesize{The comparison of the total training loss values produced by three PINN variants for Example 2.}}
\label{fig: loss2d}
\end{figure}

\begin{example}
In this example, the 2D Burgers' equations on the domain $[0, 1]^2$ with initial and boundary conditions derived from the exact solution below are solved using the proposed PINN approach. The governing equations are defined as:
\begin{flalign*}
&u_t + u u_x + v u_y = \nu(u_{xx} + u_{yy}),& \\
&v_t + u v_x + v v_y = \nu(v_{xx} + v_{yy}).&
\end{flalign*}
The system has the following exact solution \cite{caj}:
\begin{flalign*}
&u(x, y, t) = \frac{3}{4} - \frac{1}{4} \left[ 1 + \exp \left( \frac{-4x + 4y - t}{32\nu} \right) \right]^{-1},& \\
&v(x, y, t) = \frac{3}{4} + \frac{1}{4} \left[ 1 + \exp \left( \frac{-4x + 4y - t}{32\nu} \right) \right]^{-1},&
\end{flalign*}
where $x, y,$ and $t$ all belong to the interval $[0, 1]$. 
\end{example}
One of the well-known difficulties in solving CBEs with high Reynolds number $(\nu\rightarrow 0)$ is that as time advances, the solution curves steepen and develop a shock-like discontinuity. So, to compare the accuracy and effectiveness of the three PINN variants for solving the 2D CBEs with high Reynolds number $(\nu=0.0001)$, all three PINN variants utilize a consistent network architecture, the same hyperparameters, and training configurations throughout the training period as provided in Table \ref{tab:summary2}.
\begin{table}[H]
    \centering
    \begin{adjustbox}{width=\textwidth}
    \begin{tabular}{|lcccccccc|}
        \toprule
        \textbf{Method} & \textbf{Depth} & \textbf{Width} & \textbf{LR} & $N_{col}$ & $N_{ic}$ & $N_{bc}$ & \textbf{Adam Iter.} & \textbf{L-BFGS Iter.} \\
        \midrule
        LSAAL-PINN & 6 & 32 & $10^{-3}$ & 8000 & 1000 & 1000 & 4500 & 2500  \\
        AL-PINN   & 6 & 32 & $10^{-3}$ & 8000 & 1000 & 1000 & 4500 & 2500  \\
        A-PINN    & 6 & 32 & $10^{-3}$ & 8000 & 1000 & 1000 & 7000 & - \\
        \bottomrule
    \end{tabular}
    \end{adjustbox}
    \caption{Summary of hyperparameters and training configurations for the PINN variants used for Example 2.}
    \label{tab:summary2}
\end{table}
 In Figure \ref{fig: example2} (a)-(f), the predicted solution is plotted against the exact solution up to time $t=1$, and the relative $ L_2-$ Errors at four different time levels are included in the legend. From Figure \ref{fig: example2} (a)-(f), it can be seen that the LSAAL-PINN outperforms the other two PINN variants in terms of accuracy by providing about 8 times more accurate results for the component $v$ at $t=1$ than the AL-PINN and about $1385$ times more accurate than the A-PINN after $7$k runs. Moreover, one can see that the AL-PINN benefited from the hybrid optimization strategy, but lacking the self-adaptive weighting strategy, it resulted in higher errors. The A-PINN struggles most to align with the exact solution, mainly due to the lack of a dual optimization strategy and a self-adaptive weighting strategy, which results in a shock-like discontinuity in the solution of 2D CBDs. 
The spatio-temporal solution profiles of the component $u$ obtained via three PINN variants are captured in Figure \ref{fig: 2D Burger surf} with the help of surface plots up to $t=1$. As shown in Figure \ref{fig: 2D Burger surf}, the LSAAL-PINN achieves the highest alignment with the exact solution by effectively capturing the steep-gradient solution profile. Moreover, the A-PINN exhibits the lowest alignment with the exact solution and is unable to capture the shock-like discontinuity in the solution profile due to its reliance on the Adam-only optimizer. These results underscore the importance of integrating a self-adaptive weighting strategy and a dual optimization strategy to capture the sharp wave front inherent in the solution of 2D CBEs at large Reynolds numbers.
In the case of 2D CBEs, the trend of the training convergence behaviors of the three variants is depicted in Figure \ref{fig: loss2d} throughout the training process. As illustrated in the figure, the total training losses of all three variants decrease up to the magnitude of $10^{-3}$ during the initial 4500 Adam iterations. However, in the case of A-PINN, the total training losses stagnate and oscillate around the same magnitude of $10^{-3}$ after 4500 iterations, indicating no further refinement for the remainder of the training process. In contrast, the LSAAL-PINN and AL-PINN both exhibit a drastic refinement following the transition to the L-BFGS optimizer, with the LSAAL-PINN achieving a significant reduction in total training loss. These results highlight the effectiveness of the LSAAL-PINN over the other two variants for solving 2D CBEs at high Reynolds numbers.

\begin{figure}[H]
\begin{minipage}[b]{0.33\linewidth} 
\centering
\centerline{\includegraphics[width=\linewidth,height=\textheight,keepaspectratio]{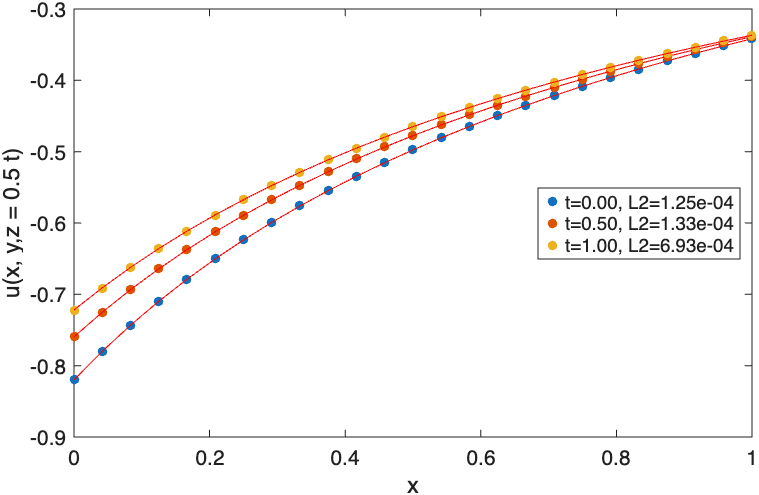}}
\centerline{\text{(a) LSAAL-PINN: u}}
\end{minipage}
\begin{minipage}[b]{0.33\linewidth} 
\centering
\centerline{\includegraphics[width=\linewidth,height=\textheight,keepaspectratio]{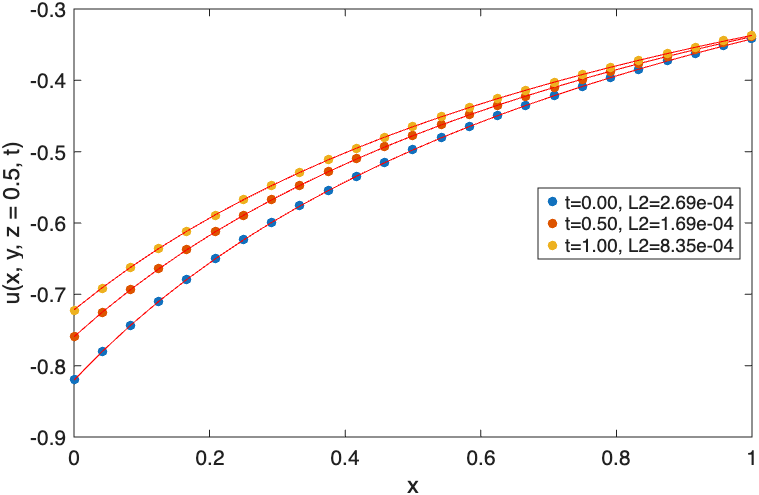}}
\centerline{\text{(b) AL-PINN: u}}
\end{minipage}
\begin{minipage}[b]{0.33\linewidth} 
\centering
\centerline{\includegraphics[width=\linewidth,height=\textheight,keepaspectratio]{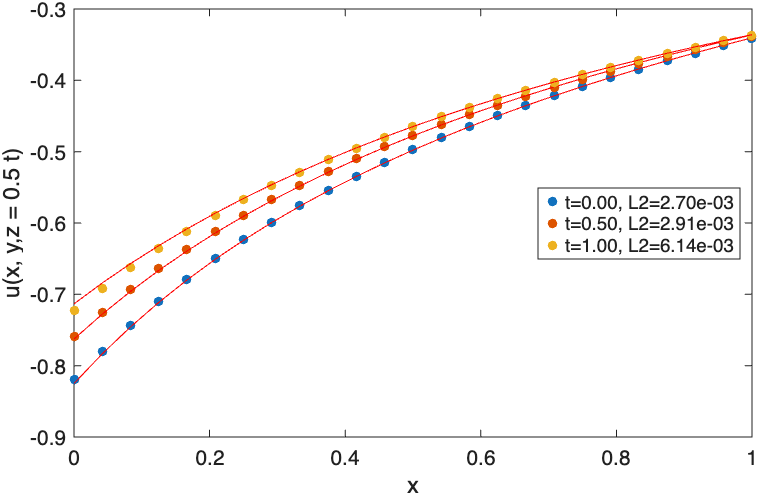}}
\centerline{\text{(c) A-PINN: u}}
\end{minipage}
\begin{minipage}[b]{0.33\linewidth} 
\centering
\centerline{\includegraphics[width=\linewidth,height=\textheight,keepaspectratio]{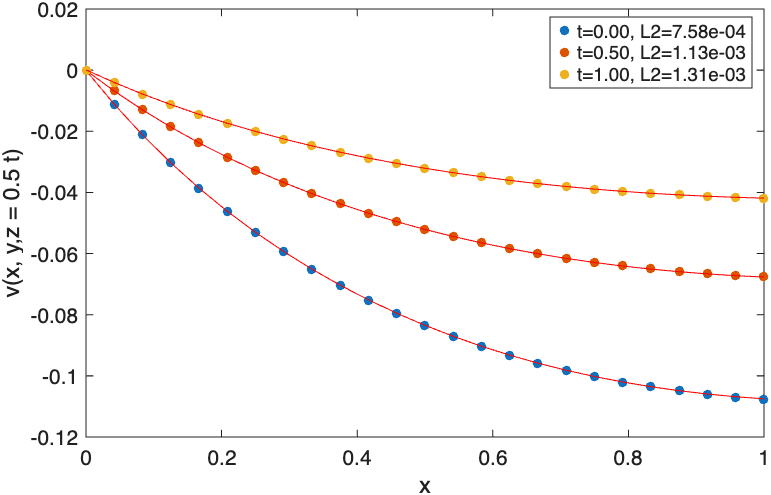}}
\centerline{\text{(d) LSAAL-PINN: v}}
\end{minipage}
\begin{minipage}[b]{0.33\linewidth} 
\centering
\centerline{\includegraphics[width=\linewidth,height=\textheight,keepaspectratio]{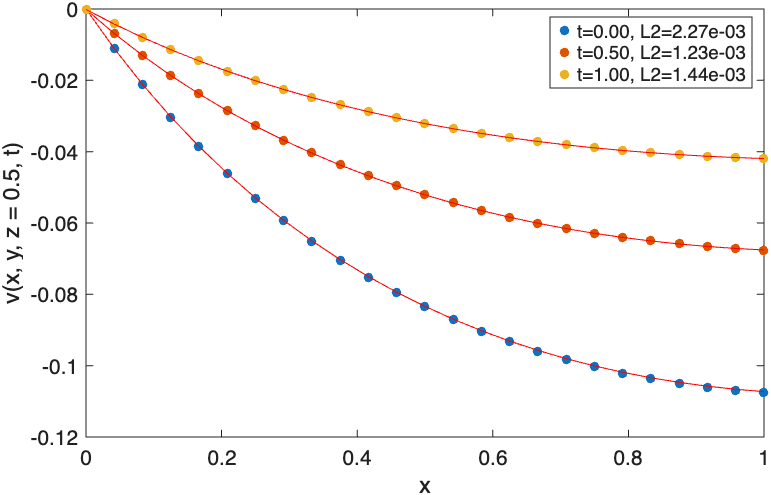}}
\centerline{\text{(e) AL-PINN: v}}
\end{minipage}
\begin{minipage}[b]{0.33\linewidth} 
\centering
\centerline{\includegraphics[width=\linewidth,height=\textheight,keepaspectratio]{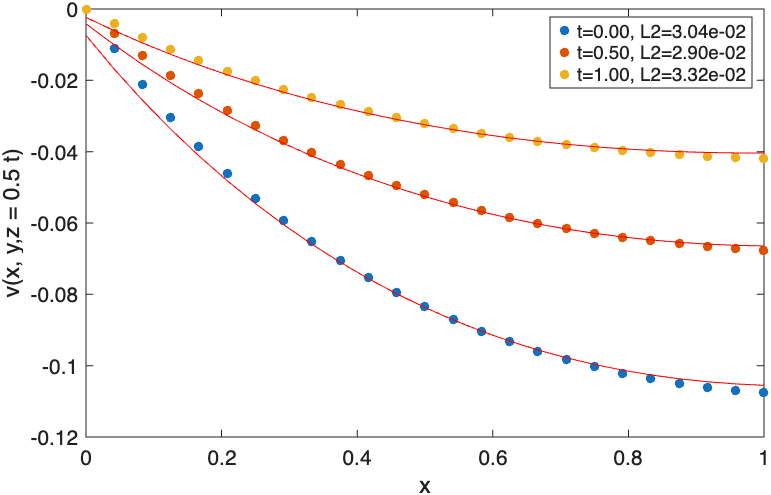}}
\centerline{\text{(f) A-PINN: v}}
\end{minipage}
\begin{minipage}[b]{0.33\linewidth} 
\centering
\centerline{\includegraphics[width=\linewidth,height=\textheight,keepaspectratio]{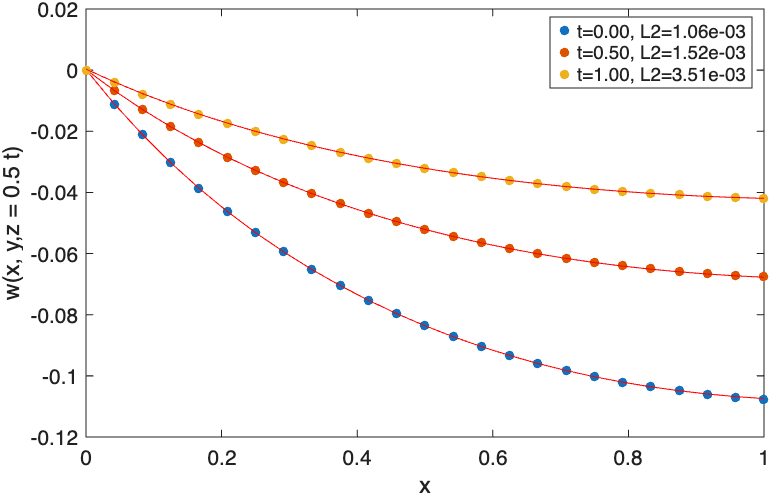}}
\centerline{\text{(g) LSAAL-PINN: w}}
\end{minipage}
\begin{minipage}[b]{0.33\linewidth} 
\centering
\centerline{\includegraphics[width=\linewidth,height=\textheight,keepaspectratio]{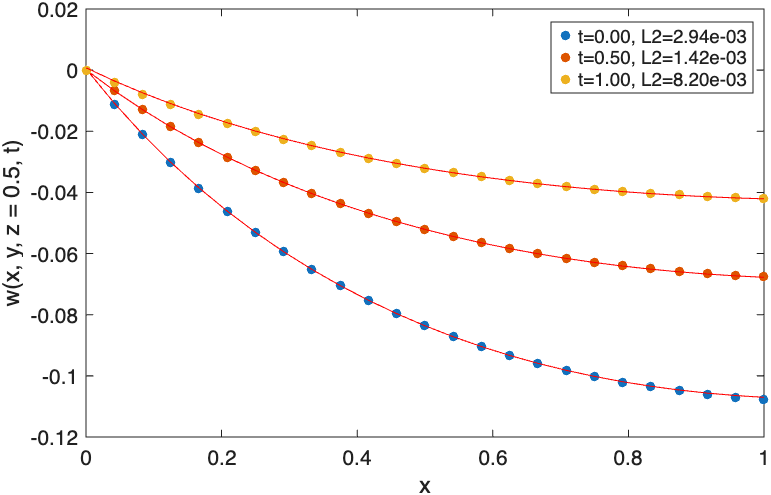}}
\centerline{\text{(h) AL-PINN: w}}
\end{minipage}
\begin{minipage}[b]{0.33\linewidth} 
\centering
\centerline{\includegraphics[width=\linewidth,height=\textheight,keepaspectratio]{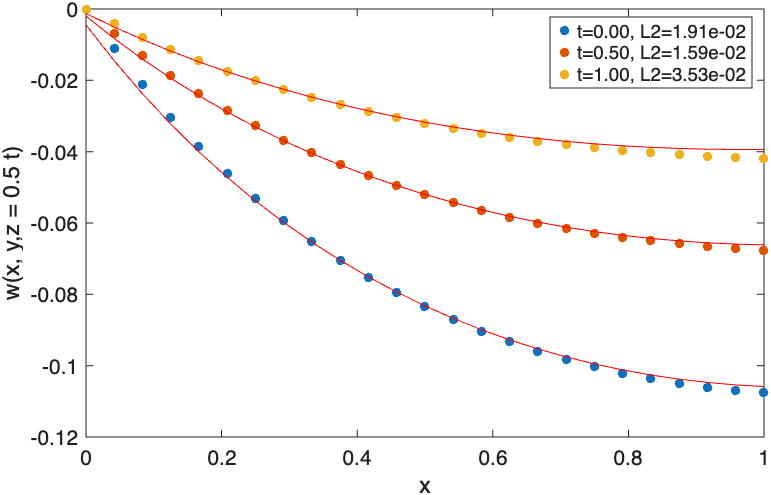}}
\centerline{\text{(i) A-PINN: w}}
\end{minipage}
 \caption{\footnotesize{Comparison of the exact solution (dots) and the predicted solution (solid lines) at various time levels for 3D CBEs. Relative $L_2-$errors for each temporal snapshot are included in the legend. From top to bottom, the first row shows the solution for $u$, the second for $v$, and the third for $w$ at $y,z=0.5$.}}
\label{fig: example3fig}
\end{figure}

\begin{figure}[H]
\begin{minipage}[b]{0.33\linewidth} 
\centering
\centerline{\includegraphics[width=\linewidth,height=\textheight,keepaspectratio]{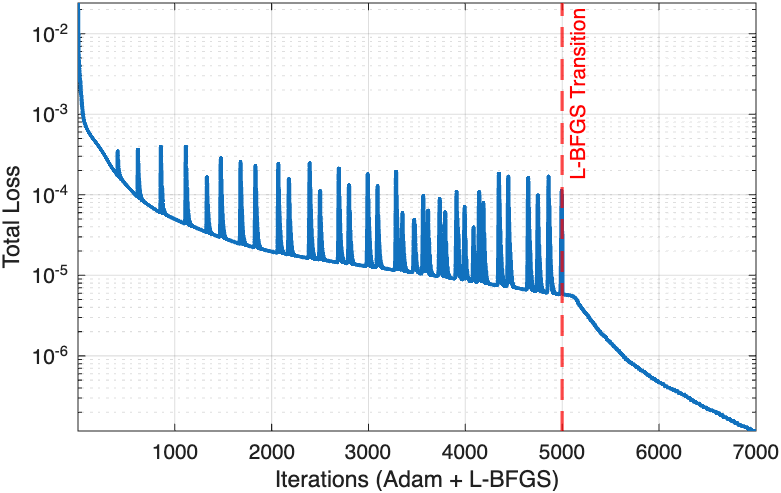}}
\centerline{\text{(a) LSAAL-PINN}}
\end{minipage}
\begin{minipage}[b]{0.33\linewidth} 
\centering
\centerline{\includegraphics[width=\linewidth,height=\textheight,keepaspectratio]{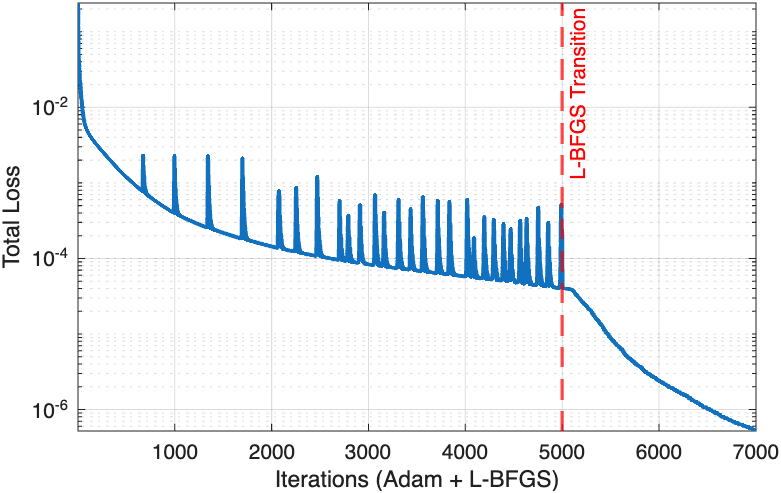}}
\centerline{\text{(b) AL-PINN}}
\end{minipage}
\begin{minipage}[b]{0.33\linewidth} 
\centering
\centerline{\includegraphics[width=\linewidth,height=\textheight,keepaspectratio]{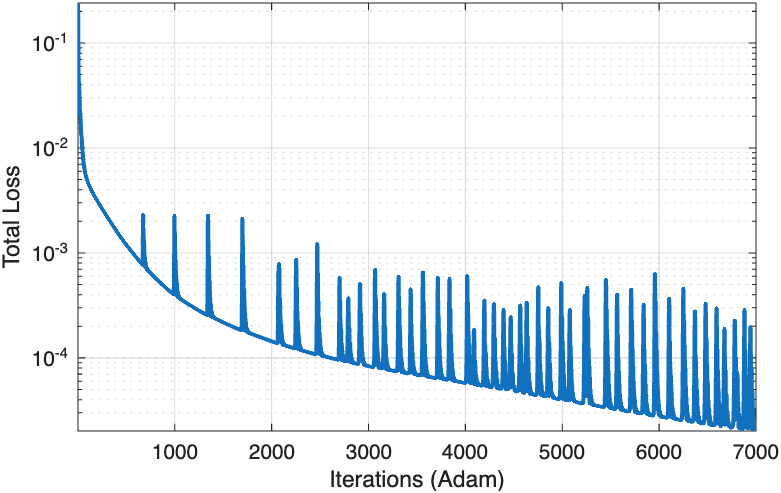}}
\centerline{\text{(c) A-PINN}}
\end{minipage}
\caption{\footnotesize{The comparison of the total training loss values produced by three PINN variants for Example 3.}}
\label{fig: loss3d}
\end{figure}

\begin{example}
For the final test problem, we considered the following 3D CBEs on domain $[-1, 1]^3$:
\begin{flalign*}
    &u_t + u u_x + v u_y + w u_z = \nu \Delta u,& \\
    &v_t + u v_x + v v_y + w v_z = \nu \Delta v,& \\
    &w_t + u w_x + v w_y + w w_z = \nu \Delta w,&
\end{flalign*}
with the initial and boundary conditions extracted from the exact solution \cite{hsm}:
\begin{flalign*}
&u(x, y, z, t) = -\frac{2}{\text{Re}} \left( \frac{\phi_x}{\phi} \right),& \\
&v(x, y, z, t) = -\frac{2}{\text{Re}} \left( \frac{\phi_y}{\phi} \right),& \\
&w(x, y, z, t) = -\frac{2}{\text{Re}} \left( \frac{\phi_z}{\phi} \right),&
\end{flalign*}
where 
\begin{flalign*}
    &\phi = 1 + x + \sin(x)\sin(y)\sin(z)e^{-t},~~\text{Re}=\frac{1}{\nu}.& 
\end{flalign*}
\end{example}
Using the information provided in the Table \ref{tab:summary3}, three PINN variants were compared to check the accuracy and effectiveness for solving 3D CBEs with $\text{Re}=3.$
\begin{table}[H]
    \centering
    \begin{adjustbox}{width=\textwidth}
    \begin{tabular}{|lcccccccc|}
        \toprule
        \textbf{Method} & \textbf{Depth} & \textbf{Width} & \textbf{LR} & $N_{col}$ & $N_{ic}$ & $N_{bc}$ & \textbf{Adam Iter.} & \textbf{L-BFGS Iter.} \\
        \midrule
        LSAAL-PINN & 6 & 32 & $10^{-3}$ & 6000 & 2500 & 2500 & 5000 & 2000  \\
        AL-PINN   & 6 & 32 & $10^{-3}$ & 6000 & 2500 & 2500 & 5000 & 2000  \\
        A-PINN    & 6 & 32 & $10^{-3}$ & 6000 & 2500 & 2500 & 7000 & - \\
        \bottomrule
    \end{tabular}
    \end{adjustbox}
    \caption{Summary of hyperparameters and training configurations for the PINN variants used for Example 3.}
    \label{tab:summary3}
\end{table}
The solution of 3D CBEs obtained via three variants was compared with the exact solution, and $L_2-$errors at various time levels were provided in the legend of the Figure \ref{fig: example3fig}(a)-(i). As we can see from the figures, the LSAAL-PINN outperforms the other two PINN variants in terms of accuracy, resulting in the lowest $ L_2-$ Errors at various time levels throughout the training process. These results justify the effectiveness of the proposed framework for predicting the solution profile of the 3D DBEs more accurately than AL-PINN and A-PINN. Consequently, highlight the importance of hybrid optimizers and self-adaptive weighting strategies for solving 3D CBEs.
In the case of 3D CBEs, the performance differences across the three PINN variants are depicted in Figure \ref{fig: loss3d} through the loss plots. As illustrated in the figure, the total training loss of A-PINN at $7$k runs reaches around the magnitude of $10^{-5}$, which is even larger than the loss value achieved by the LSAAL-PINN at $5$k runs. While the AL-PINN was able to refine the loss value by implementing L-BFGS after $5$k iterations, its loss value after $7$k iterations remains higher than that of the LSAAL-PINN. Overall, the LSAAL-PINN achieves a significant reduction in total training loss. These results highlight the effectiveness of the LSAAL-PINN over the other two variants for solving 3D CBEs.

\begin{figure}[H]
\begin{minipage}[b]{0.33\linewidth} 
\centering
\centerline{\includegraphics[width=\linewidth,height=\textheight,keepaspectratio]{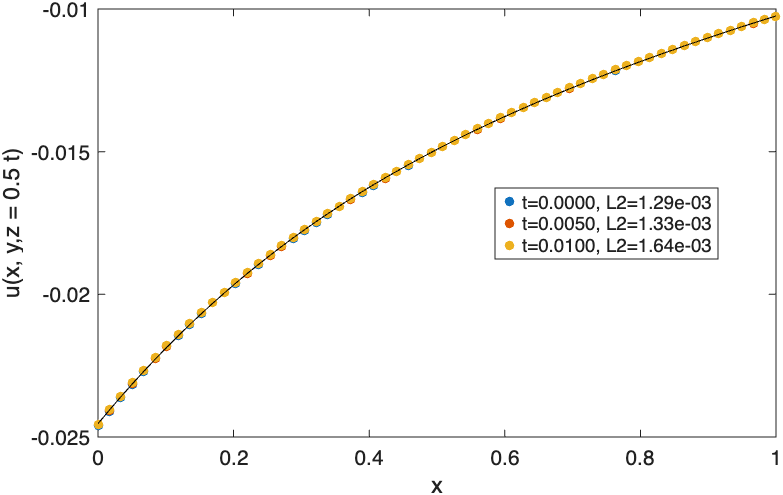}}
\end{minipage}
\begin{minipage}[b]{0.33\linewidth} 
\centering
\centerline{\includegraphics[width=\linewidth,height=\textheight,keepaspectratio]{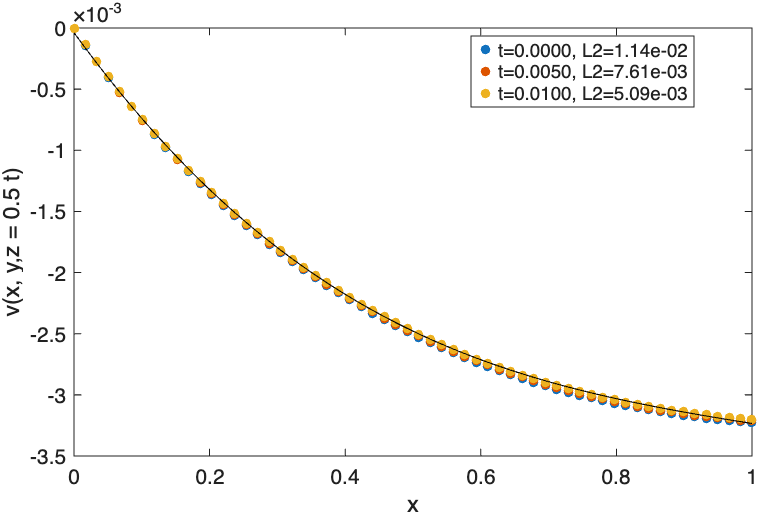}}
\end{minipage}
\begin{minipage}[b]{0.33\linewidth} 
\centering
\centerline{\includegraphics[width=\linewidth,height=\textheight,keepaspectratio]{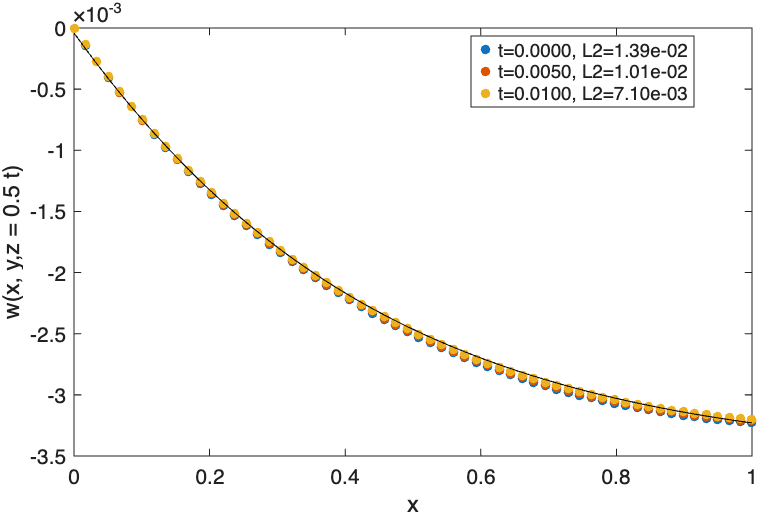}}
\end{minipage}
\caption{\footnotesize{Comparison of the exact solution (dots) and the predicted solution (solid lines) obtained via LSAAL-PINN at various time levels for 3D CBEs with $\text{Re = 100}$. Relative $L_2-$errors for each component are included in the legend.}}
\label{fig: example3l2}
\end{figure}

To check the performance of the proposed LSAAL-PINN for solving 3D CBEs with high Reynolds number $\text{Re = 100}$, a numerical experiment ran on Example 3 with $N_{\text{col}}=8500,N_{\text{ic}}=1500, N_{\text{bc}}=1500,$ Adam-Iter $=4500$, L-BFGS Iter $=2500$, same network  architecture and learning rate provided in Table \ref{tab:summary3}. From the Figure \ref{fig: example3l2}, it can be seen that the LSAAL-PINN predicts the solution of 3D CBEs with high Reynolds number reasonably well and aligns with the exact solution accurately, providing the smaller $L_2-$errors in the magnitude of $10^{-3}.$

\begin{figure}[H]
\begin{minipage}[b]{0.33\linewidth} 
\centering
\centerline{\includegraphics[width=\linewidth,height=\textheight,keepaspectratio]{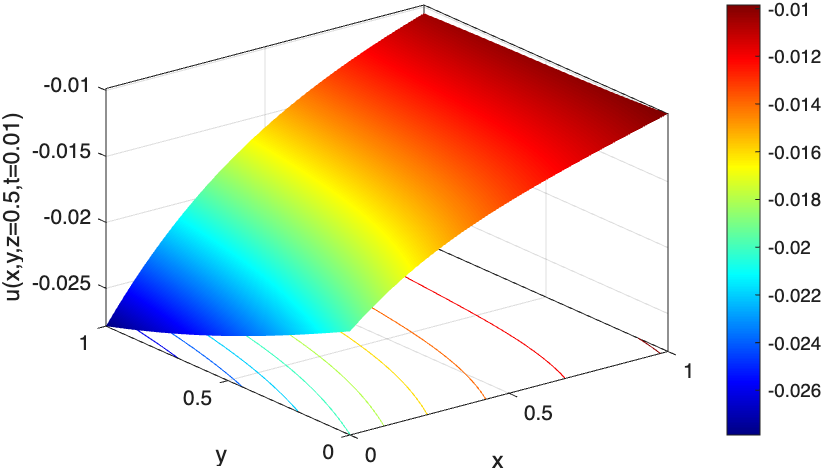}}
\centerline{\text{(a) LSAAL-PINN: u}}
\end{minipage}
\begin{minipage}[b]{0.33\linewidth} 
\centering
\centerline{\includegraphics[width=\linewidth,height=\textheight,keepaspectratio]{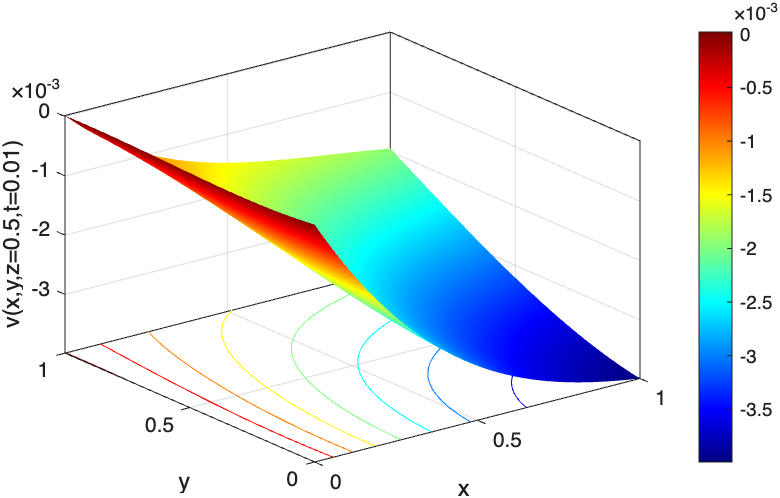}}
\centerline{\text{(b) LSAAL-PINN: v}}
\end{minipage}
\begin{minipage}[b]{0.33\linewidth} 
\centering
\centerline{\includegraphics[width=\linewidth,height=\textheight,keepaspectratio]{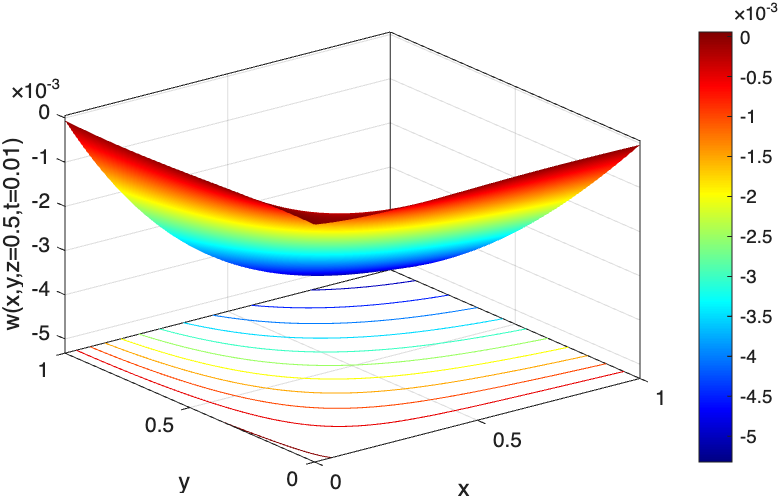}}
\centerline{\text{(c) LSAAL-PINN: w}}
\end{minipage}
\begin{minipage}[b]{0.33\linewidth} 
\centering
\centerline{\includegraphics[width=\linewidth,height=\textheight,keepaspectratio]{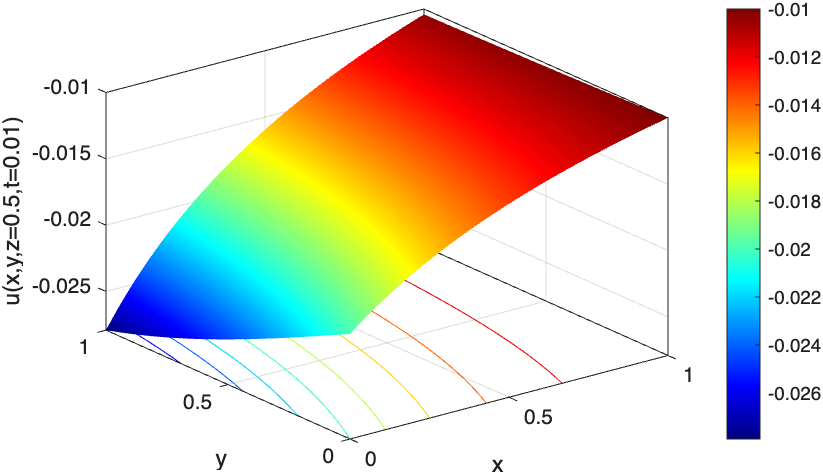}}
\centerline{\text{(d) Exact $u$}}
\end{minipage}
\begin{minipage}[b]{0.33\linewidth} 
\centering
\centerline{\includegraphics[width=\linewidth,height=\textheight,keepaspectratio]{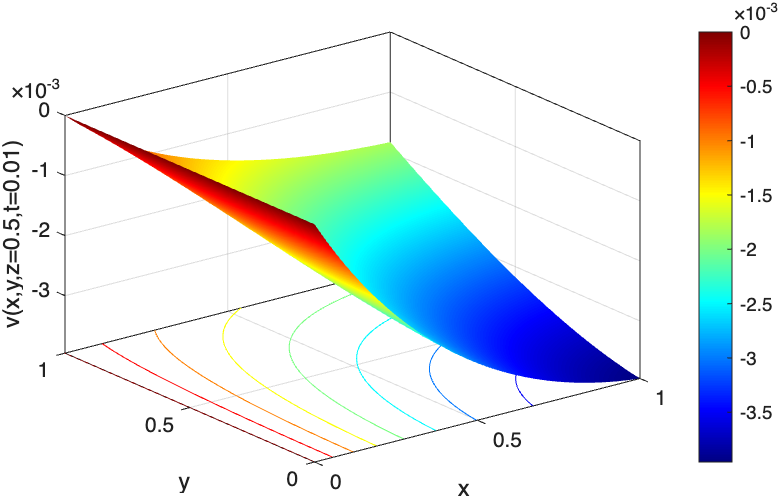}}
\centerline{\text{(e) Exact $v$}}
\end{minipage}
\begin{minipage}[b]{0.33\linewidth} 
\centering
\centerline{\includegraphics[width=\linewidth,height=\textheight,keepaspectratio]{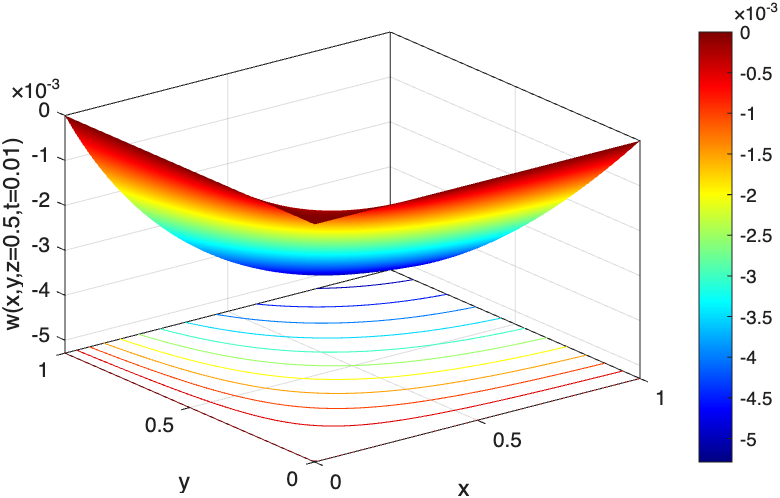}}
\centerline{\text{(f) Exact $w$}}
\end{minipage}
\caption{\footnotesize{The exact and LSAAL-PINN predicted solution of 3D CBEs with $\text{Re}=100$} at $t=0.01.$}
\label{fig: example3surf}

\end{figure}

In the Figure \ref{fig: example3surf}(a)-(f), the spatio-temporal solution profile of the 3D CBEs with high Reynolds number at $t=0.01$ obtained via the LSAAL-PINN is compared with the exact solution. As we can see from the figure, the predicted solutions are more closely aligned with the exact solutions, justifying the effectiveness of hybrid optimization and self-adaptive weighting strategies in solving 3D CBEs.

\section{Conclusion}
This study proposed a layerwise self-adaptive PINN framework with hybrid optimization strategies and compared it with a standard PINN with and without hybrid optimization strategies for solving MCBEs with high Reynolds numbers. Numerical results across 1D, 2D, and 3D domains confirmed that the LSAAL-PINN framework is more effective, accurate, and robust than the standard PINN with and without dual-phase optimization for simulating the spatio-temporal solution profile of CBEs with vanishingly small viscosity. In particular, the proposed framework is more effective in capturing the development of sharp shock fronts as time evolves in the solution of Burgers' equation. These results highlight that the integration of a layerwise self-adaptive weighting strategy and the dual-phase optimization approach is sufficient to capture the steep gradient solution profile without the need for a computationally expensive pointwise weighting strategy. Although this work is limited to the MCBEs, future research may extend the implementation of the proposed framework to the more challenging time-dependent PDEs easily.

\bibliographystyle{unsrt}

\end{document}